\documentclass{elsart}
\usepackage{amsmath,amsbsy,amsfonts}
\usepackage{natbib}
\usepackage{epsfig}

\def\mydate{\number\year.\ifnum\month<10{0}\fi\number\month.\ifnum\day<10{0}\fi\number\day}
\renewcommand{\theequation}{\thesection.\arabic{equation}}
\let\ssection=\section
\renewcommand{\section}{\setcounter{equation}{0}\ssection}

\newcommand{\dy}{\partial}
\newcommand{\ddt}[1]{\frac{{\rm d}{#1}}{{\rm d}{t}}}
\newcommand{\ddy}[2]{\frac{\partial{#1}}{\partial{#2}}}

\newcommand{\df}{{\bf d}}

\newcommand{\lapl}{\Delta}
\newcommand{\eps}{\varepsilon}
\newcommand{\vfi}{\varphi}
\newcommand{\RReal}{\mathbb{R}}
\newcommand{\Comp}{\mathbb{C}}
\newcommand{\Zahl}{\mathbb{Z}}
\newcommand{\Oh}{{O}}

\newcommand{\gb}{\boldsymbol{\nabla}}
\newcommand{\sgb}{\boldsymbol{\nabla}^\perp}

\newcommand{\xb}{\boldsymbol{x}}
\newcommand{\Xb}{\boldsymbol{X}}
\newcommand{\vb}{\boldsymbol{v}}
\newcommand{\dxdy}{\;{\rm d}x\,{\rm d}y}

\newcommand{\sfrac}[2]{{\textstyle\frac{#1}{#2}}}

\newcommand{\ds}{\;{\rm d}s}

\newcommand{\psiL}{\psi_\Lambda}
\newcommand{\wL}{\omega_\Lambda}

\newcommand{\Lz}{{0}}

\newcommand{\geop}{\Delta\theta_{\rm geo}}
\newcommand{\dynp}{\Delta\theta_{\rm dyn}}

\newcommand{\cpoi}{c_{{\rm poi}}}

\newcommand{\PpsiL}{{\sf P}_{\Lambda}}
\newcommand{\Ppsiz}{{\sf P}_{0}}

\newcommand{\omegaz}{\omega^{(0)}}
\newcommand{\psiz}{\psi^{(0)}}
\newcommand{\psizb}{\hat\psi_\Lambda}
\newcommand{\omegao}{\omega^{(1)}}
\newcommand{\psio}{\psi^{(1)}}
\newcommand{\Psio}{\Psi^{(1)}}
\newcommand{\Psioh}{\hat\Psi^{(1)}}
\newcommand{\Psiob}{\bar\Psi^{(1)}}
\newcommand{\Lambdadot}{\dot\Lambda}
\newcommand{\Lambdaz}{0}

\newcommand{\dLwdL}[2]{\df\Lambda_{#1}\wedge\df\Lambda_{#2}}

\newcommand{\gL}{g_\Lambda}
\newcommand{\GL}{G_\Lambda}
\newcommand{\FL}{F_\Lambda}

\newcommand{\Lspace}{\mathcal{L}}
\newcommand{\Area}{\mathcal{A}}

\newcommand{\CFun}{\mathcal{C}}

\newcommand{\phiform}{\Phi}
\newcommand{\phinat}{\Phi^\star}
\newcommand{\phinath}{\hat{\Phi}^\star}
\newcommand{\phinatb}{\bar{\Phi}^\star}
\newcommand{\phiformb}{\hat{\Phi}}

\newcommand{\Psibc}{\boldsymbol{\beta}}
\newcommand{\Psib}[1]{\bar{\Psi}^{(#1)}}

\def\d2x{\;\mathrm{d}^2\xb}
\newcommand{\sig}{\sigma}
\newcommand{\ex}{\mathrm{e}}
\renewcommand{\i}{\mathrm{i}}
\def\Pinat{\Pi^\star}

\begin{document}

\begin{frontmatter}

\title{Two-dimensional Euler flows in slowly deforming domains}

\author[Vanneste]{J.~Vanneste\corauthref{cor}}
\corauth[cor]{Corresponding author.}
\address[Vanneste]{School of Mathematics and Maxwell Institute for Mathematical Sciences, University of Edinburgh, Edinburgh~EH9~3JZ, UK}
\ead{j.vanneste@ed.ac.uk}

\author[Wirosoetisno]{D.~Wirosoetisno}
\address[Wirosoetisno]{Department of Mathematical Sciences, University of Durham\\ Durham~DH1~3LE, UK}
\ead{djoko.wirosoetisno@durham.ac.uk}

\begin{abstract}
We consider the evolution of an incompressible two-dimensional perfect fluid
as the boundary of its domain is deformed in a prescribed fashion.
The flow is taken to be initially steady, and the boundary deformation
is assumed to be slow compared to the fluid motion.
The Eulerian flow is found to remain approximately steady throughout
the evolution.
At leading order, the velocity field depends instantaneously on the shape
of the domain boundary, and it is determined by the steadiness and
vorticity-preservation conditions.
This is made explicit by reformulating the problem in terms of an
area-preserving diffeomorphism $\gL$ which transports the vorticity.
The first-order correction to the velocity field is linear in the
boundary velocity, and we show how it can be computed from the time-derivative
of $\gL$. 

The evolution of the Lagrangian position of fluid particles is also examined.
Thanks to vorticity conservation, this position can be specified by an
angle-like coordinate along vorticity contours.
An evolution equation for this angle is derived, and the net change in
angle resulting from a cyclic deformation of the domain boundary is
calculated.
This includes a geometric contribution which can be expressed as the
integral of a certain curvature over the interior of the circuit that
is traced by the parameters defining the deforming boundary.

A perturbation approach using Lie series is developed for the
computation of both the Eulerian flow and geometric angle for small
deformations of the boundary.
Explicit results are presented for the evolution of nearly axisymmetric
flows in slightly deformed discs.
\end{abstract}


\end{frontmatter}

% ===========================================================================

\section{Introduction}

This paper examines the dynamics of a two-dimensional (2d) fluid inside
a container whose boundary is deformed slowly.
The fluid is assumed to be perfect and incompressible; consistent with
the latter assumption, the area of the container is constant.
Beyond potential applications such as the control of fluid flows, we use the
problem as a paradigm for the study of Hamiltonian fluid models depending
on slowly varying parameters.
This is an obvious first step:
the 2d Euler equations governing incompressible perfect fluids are
indeed Hamiltonian \cite[e.g.,][]{pjm:98,salmon:88}, and imposing
boundary deformations is arguably the most natural way of introducing
a parameter dependence.
As is well-known in finite dimensions, Hamiltonian systems are
strongly constrained; as a result, slow changes of parameters lead
to a remarkable behaviour encapsulated in the theory of adiabatic
invariance \cite[cf., e.g.,][]{arnold:mmcm,landau-lifshitz:m} and
geometric angles \cite[]{hannay:85,berry:85}.
In 2d Euler, the material invariance of vorticity \cite[e.g.,][]{saffman:vd}
similarly imposes a strong contraint on the system, which we exploit
extensively to derive what can be interpreted as fluid-dynamical versions
of adiabatic invariance and geometric angle.

The problem we consider here is rather involved in its full generality.
To make progress, we make a number of assumptions and consider the
following scenario.
At an initial time, a steady flow is given in some simply-connected
domain $D_\Lambdaz$.
The streamlines have the simplest topology, that of nested closed
curves, and the flow is Arnold stable (see section \ref{s:eu0} below).
We then assume that this continues to hold throughout the evolution
as the domain is being deformed.
With these hypotheses, we use an asymptotic approach, based on the
separation between the timescales of the boundary deformation
and that of the flow, to answer two questions:
(i)~what is the leading-order approximation to the (Eulerian) flow at any time;
and (ii)~what is the (Lagrangian) position of the fluid particles? 

\nocite{dw-jv:ddef}
The first question is answered by showing that the leading-order flow
is steady at all times.
This makes it possible to rephrase the problem in terms of an
area-preserving diffeomorphism $\gL$ which maps the vorticity in
the initial domain to the vorticity in the deformed domain.
The uniqueness of $\gL$ up to displacements along lines of constant
vorticity, established in Wirosoetisno and Vanneste (2005; henceforth WV)
and revisited here, shows that the leading-order velocity field is
completely determined by the instantaneous shape of the boundary
and is independent of the history of past shapes.
This, of course, is analogous to the adiabatic invariance of the action variables for
finite-dimensional Hamiltonian systems with slowly varying parameters. 

To answer the second question, we need to compute the first-order
correction to the approximate velocity field obtained in~(i).
This is because the evolution equation for the particle position needs to be
integrated over the long time scales required to achieve order-one boundary
deformations.
The first-order correction to the velocity field is linear in the boundary
velocity, and it can be derived from $\mathrm{d}\gL/\mathrm{d}t$
by solving a pseudodifferential equation.
Once this is done, the particle-position problem reduces to the
solution of (independent) one-degree-of-freedom Hamiltonian
systems with slowly varying parameters. 
Since particles remain on vorticity contours (which in effect
are contours of constant action), only the position along each
contour, regarded as an angle variable, needs to be determined.
The value of this angle is found to depend on the history of the
boundary shape. It includes a geometric contribution, similar to
the Hannay--Berry angle, which possesses a nice interpretation
\cite[]{hannay:85,berry:85,montgomery:88,shapere-wilczek:89,marsden-montgomery-ratiu:90}.
We note that the geometric angle has been studied in fluid dynamics
by \cite{shashikanth-newton:97,shashikanth-newton:99} who
considered point-vortex solutions of  the 2d Euler equation, and
by \cite{shapere-wilczek:jfm:89} for Stokes flow; here it appears
in the context of smooth inviscid flows.

The determination of the leading-order Eulerian flow from the steadiness
and vorticity-preservation conditions was treated in WV,
where conditions for the existence of $\gL$ and the uniqueness of the
resulting velocity field were given in appropriate function spaces
for sufficiently small boundary deformations.
In the present article we adopt a more
informal approach to treat both the Eulerian and Lagrangian problems
under a slightly different set of hypotheses;
rigorous proof of the adiabatic invariance will be the subject of
a future work.
It proves convenient to express our derivation in the language of
differential forms in the space of the parameters defining the boundary shape.
This makes explicit the linear dependence of several important
quantities on the boundary velocity, and it gives a natural description
of the geometric angle in terms of a curvature form in the parameter space.
We use this language mainly as a notational tool, but it is clear that
a more abstract geometric interpretation of the results could be given.
This is discussed at the end of the paper.

In the following section, we present a short description of the
2d Euler equation in a deforming domain in order to fix the notation,
and we consider the behaviour of the leading-order Eulerian
flow.
Next, in \S\ref{s:eu1} we compute the first-order correction to the
Eulerian flow, which depends only on the instantaneous shape of the
boundary and its velocity.
Using these results, we study the Lagrangian flow in \S\ref{s:lag}, where
the geometric angle of the particle position is derived.
In these sections, we consider general domains and arbitrary
boundary deformations, requiring only that the boundary deformation
be slow.
The results are given as solutions of partial (pseudo)differential
equations, which in general will have to be solved numerically.
In \S\ref{s:small} we develop a perturbative approach for the solution
of these equations, based on the assumption of small (total) boundary
deformation.
We carry out the calculation to second order, but the Lie series
formulation that we employ is well suited for systematic extensions
to higher orders.
An application to nearly axisymmetric flows in a slightly deformed
disc is presented in \S\ref{s:axi},
followed by a Discussion.
Technical details are relegated to the Appendices.

% \medskip\hbox to\hsize{\qquad\hrulefill\qquad}\medskip

\begin{figure}\label{f:LDL}
\begin{center}
\epsfig{file=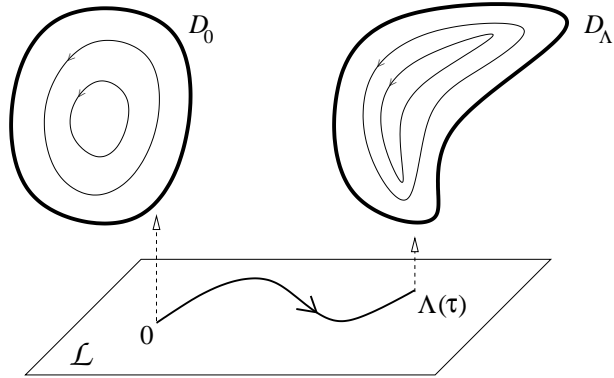,height=5cm}
\caption{Parameterization of the shape of the domain $D_\Lambda$ by $\Lambda$:
as the parameter $\Lambda$ moves from $\Lambda(0)=0$ to $\Lambda(\eps t)$
in $\Lspace$, the fluid's domain changes its shape from $D_0$ to
$D_\Lambda$, inducing a change in the leading-order Eulerian flow
whose streamlines are indicated.}
\end{center}
\end{figure}

% ===========================================================================

\section{Eulerian Flow: Adiabatic Invariance}\label{s:eu0}

We begin by studying the behaviour of the Eulerian flow.

% ---------------------------------------------------------------------------

\subsection{Formulation}

Let $D_{\Lambda(\eps t)}\subset\RReal^2$ be a simply-connected, bounded
and smooth domain which is slowly evolving in time $t$ in a prescribed
fashion while keeping its area fixed. Here, $\Lambda$ denotes the set of
parameters defining the shape of the boundary, and the slowness of their
time dependence is made explicit by the introduction of the asymptotic
parameter $0<\eps\ll 1$. 
Denoting the (generally infinite-dimensional) space in which $\Lambda$
lives by $\Lspace$, we can think of the evolution of the domain shape
as the tracing of a curve $\Lambda(\eps t) \subset \Lspace$;
see Figure~\ref{f:LDL}.
In this article, our concern is the behaviour of the flow as $\eps\to0$,
that is, in the limit of slow boundary deformation, over $\Oh(\eps^{-1})$
time scales so that $\Oh(1)$ deformations are achieved.
We make a blanket assumption that all functions are sufficiently
smooth for our purposes, and we denote by $\CFun(D_\Lambda)$ the
space of smooth real-valued functions in $D_\Lambda$.

We can describe the evolution of a perfect fluid flow in
$D_{\Lambda(\eps t)}$ using the vorticity--streamfunction formulation
\begin{align}
   &\dy_t\omega + [\psi,\omega] = 0, \label{q:dwdt}\\
   &\omega = \lapl\psi.
\end{align}
The velocity is given by $(u,v)=\sgb\psi:=(-\dy_y\psi,\dy_x\psi)$,
$\omega=\sgb\cdot\vb:=\dy_x v-\dy_y u$ is the vorticity, and
$[f,g] := \sgb f \cdot \gb g = \dy_x f\,\dy_y g - \dy_x g\,\dy_y f$ is
the Jacobian.
Steady flows satisfy $[\psi,\omega]=0$. 

A convenient way of of defining the domain boundary $\dy D_{\Lambda(\eps t)}$
is as the level set $B(x,y;\Lambda(\eps t))=0$ of some prescribed function $B$.
Since $\dy D_{\Lambda(\eps t)}$ is a material curve, 
\begin{equation}\label{q:dBdt}
   \dy_t B + [\psi,B] = 0
   \qquad\textrm{for } B(x,y;\Lambda(\eps t))=0.
\end{equation}
Assuming that $\gb B\ne0$ on $\dy D_\Lambda$, this can be inverted to give
the boundary condition
\begin{equation}\label{q:bc0}
   \psi = \eps\, b(\eps t) \qquad\textrm{on }\dy D_{\Lambda(\eps t)}.
\end{equation}
Since $\psi$ is determined only up to an additive constant, we set
\begin{equation}
   \oint_{\dy D_{\Lambda(\eps t)}} b(\eps t) \;{\rm d}l = 0.
\end{equation}
It is clear from \eqref{q:dBdt}--\eqref{q:bc0} that $\psi$ is proportional to
${\rm d}\Lambda/{\rm d}t$ on $\dy D_\Lambda$.

Exploiting the smallness of $\eps$, we expand the vorticity
and streamfunction in $\eps$,
\begin{equation}\label{q:wpeps}
   \omega = \omegaz + \eps\,\omegao + \eps^2\omega^{(2)} + \cdots
   \quad\textrm{and}\quad
   \psi = \psiz + \eps\,\psio + \eps^2\psi^{(2)} + \cdots.
\end{equation}
Our aim in this section is to compute the leading-order flow $\psiz$,
given its initial value and the boundary deformation $b(\eps t)$,
and to show that it depends only on the boundary shape $\Lambda$ and
not on its time history.

First we note that the boundary conditions \eqref{q:bc0} imply that on $\dy D_\Lambda$,
\begin{equation}\label{q:psibc}
   \psio = b
   \qquad\textrm{and}\qquad
   \psi^{(n)} = 0\textrm{  for }n=0,2,3,\cdots.
\end{equation}
Since the total vorticity $\omega$ is advected by the flow,
the boundary $\dy D_\Lambda$ is a vorticity contour and thus we can take
$\omega=\omegaz$ and $\omega^{(n)}=0$ for $n=1,2,\cdots$ there.

Substituting \eqref{q:wpeps} into \eqref{q:dwdt}, we find
\begin{equation}\label{q:eomt}
   [\psiz,\omegaz]
   + \dy_t\omegaz + \eps\,[\psio,\omegaz] + \eps\,[\psiz,\omegao]
   + \eps\,\dy_t\omegao + \Oh(\eps^2) = 0.
\end{equation}
If the fluid flow is stable in the absence of boundary deformation,
we expect that the flow will evolve only slowly when the boundary
is deforming.
We therefore introduce the slow time
\begin{equation}
\tau=\eps t,
\end{equation}
in terms of which 
\eqref{q:eomt} becomes
\begin{equation}
   [\psiz,\omegaz]
   + \eps \dy_\tau\omegaz + \eps\,[\psio,\omegaz] + \eps\,[\psiz,\omegao]
   + \Oh(\eps^2) = 0.
\end{equation}
At leading order we obtain
\begin{equation}\label{q:psi0w0}
   [\psiz,\omegaz] = 0.
\end{equation}
Taking into account the fact that $\psiz=0$ on $\dy D$, we find that
the leading-order flow $\psiz$ is {\em instantaneously steady\/}.
The relation \eqref{q:psi0w0} implies that there exists a scalar function
$G$ relating $\omegaz$ and $\psiz$,
\begin{equation}\label{q:Psi0}
   \psiz = G(\omegaz;\tau).
\end{equation}
As noted, the function $G$ depends on the slow time $\tau$,
regarded here as a parameter for reasons which will be apparent later.
We define $F$ as the inverse of $G$:
$G(F(u;\tau);\tau)=u$ for every $u$ and $\tau$.
With an abuse of notation, we will often write
$G'$ for $G'\circ\omegaz=\gb\psiz/\gb\omegaz$ and
$F'$ for $F'\circ\psiz=\gb\omegaz/\gb\psiz$;
what is meant will be clear from the context.

At $\Oh(\eps)$ we have
\begin{equation}\label{q:dw0dt1}
   \dy_\tau\omegaz + [\psio,\omegaz] + [\psiz,\omegao] = 0.
\end{equation}
Using \eqref{q:Psi0}, the second term can be written as
\begin{equation}
   [\psiz,\omegao] = G' [\omegaz,\lapl\psio]
   = [\omegaz,G'\lapl\psio],
\end{equation}
with which \eqref{q:dw0dt1} becomes
\begin{align}
  &\dy_\tau\omegaz + [\phi,\omegaz] = 0,  \label{q:dw0dt}\\
  &\phi = [1 - G'\lapl]\,\psio.   \label{q:phidef}
\end{align}
These two equations imply that the leading-order vorticity $\omegaz$
is {\em rearranged\/} by a divergence-free velocity field $\sgb\phi$
with $\phi$ related to the first-order streamfunction $\psio$
by \eqref{q:phidef}.

% ---------------------------------------------------------------------------

\subsection{Determination of the Eulerian flow} 

We now show how the leading order flow $\psiz$ can be determined from
\eqref{q:Psi0} and the fact that $\omegaz(t)$ is a rearrangement of its
initial value $\omegaz(0)$.
We make the following two assumptions on $\psiz$:

\medskip\noindent{\bf H1.}
{\sl The leading-order streamfunction $\psiz$ is such that it has
only one critical point in $D_\Lambda$ (which is necessarily elliptic)
and is nonlinearly stable in the sense of Arnold.\/}

\smallskip\noindent
We recall that Arnold stability [cf.\ \cite{holm-al:stab:85}] requires that
the steady streamfunction $\psiz$ satisfies either
(i)~$0<c_1\le G'\le c_2<\infty$, or
(ii)~$0<1/\cpoi<c_1\le -G'\le c_2<\infty$.
In the second condition, $\cpoi$ is the Poincar{\'e} constant for the
domain $D_\Lambda$, namely the smallest eigenvalue $\mu$ of the problem
\begin{equation}
   (\Delta+\mu)\, u = 0 \textrm{ in } D_\Lambda
   \qquad\textrm{with}\qquad
   u = 0 \textrm{ on }\dy D_\Lambda.
\end{equation}
These conditions ensure that the steady flow is either a minimum or
a maximum of the energy for fixed vorticity distribution.
Note that H1 implies that 
\begin{equation}
-\cpoi < F' < \infty,
\end{equation}
a condition which will be useful below.
The assumption H1 is stronger than that made in WV but it considerably
simplifies the solution of \eqref{q:pde0}--\eqref{q:chibc} below.

For the second assumption, we need a little more notation.
Let $s$ denote a variable conjugate to $\psi^{(0)}$ in $D_\Lambda$, satisfying 
$[\psi^{(0)},s]=1$. Denoting the differential arclength along the
curve $\psi^{(0)}=\textrm{const}$ by ${\rm d}l$, we have
${\rm d}s={\rm d}l/|\gb\psi^{0}|$.
We then assume:

\medskip\noindent{\bf H2.}
{\sl There exists a $c_\psi>0$ such that, for all values of $c$
assumed by $\psi^{(0)}$,
\begin{equation} \label{q:h2}
   \oint_{\psi^{(0)}=c} \ds \le \frac{1}{c_\psi}.
\end{equation}}

\medskip\noindent
In the context adiabatic invariance,
this condition is natural: the left-hand side of \eqref{q:h2} gives the
period of rotation of fluid parcels along the streamline $\psi^{(0)}=c$;
its boundedness ensures that a time-scale separation between this period
and the time scale of the boundary deformation exists for sufficiently
small $\eps$.
As noted in WV, H2 holds if $\wL\ne0$ at the fixed point of $\psiL$.

\medskip
For concreteness, we henceforth assume that, at $t=0$, our domain is
parameterised by $\Lambda_0$ and we choose our coordinates in $\Lspace$
such that $\Lambda_0=\Lambdaz$.
Furthermore, we fix in $D_{\Lambda_0}\equiv D_\Lambdaz$ a steady
leading-order flow $\psi^{(0)}(\xb,0)=\psi_\Lambdaz(\xb)$ satisfying H1--H2.
Considering only the leading-order flow $\psi^{(0)}(\xb,t)$ for the moment,
we then claim that, assuming H1--H2:

\medskip\noindent{\bf P1.}
{\sl
The flow $\psi^{(0)}$ is uniquely determined by
(i)~the shape of the deformed domain $D_\Lambda$,
(ii)~the steadiness condition \eqref{q:Psi0}, and
(iii)~the fact that the vorticity $\omega^{(0)}=\Delta \psi^{(0)}$
is obtained by rearrangement of the initial vorticity
$\omega_{\Lambdaz} = \Delta\psi_\Lambdaz$.} 

\smallskip\noindent
As a result, the leading-order flow at a fixed time $t$ depends only
on the shape of the deformed domain at $t$ (parameterized by
$\Lambda(\eps t)$), and not on the history of shapes at intermediate
times (parameterized by the path $\Lambda(\tau),\,0<\tau<\eps t$).
One may draw an analogy with adiabatic invariance in finite-dimensional
Hamiltonian systems with slowly-varying parameters: here the amplitude
is completely determined by the instantaneous value
of the parameter, not by its time history.

To emphasize the fact that the leading-order flow depends on $\Lambda$
instantaneously, we introduce the notation
\begin{equation}
   \psi^{(0)} = \psiL \quad \textrm{and} \quad \omega^{(0)} = \wL
\end{equation}
for the leading-order streamfunction and vorticity.
We also write $G(\cdot;\tau)=:G_{\Lambda(\tau)}(\cdot)$ and
$F(\cdot;\tau)=:F_{\Lambda(\tau)}(\cdot)$.
These define the scalar functions $\GL$ and $\FL$, both of which have
$\Lambda$ as a parameter.
When (and only when) no confusion may arise, we will often write $\GL\circ\wL$
as $\GL$, $\GL'\circ\wL$ as $\GL'$, and similarly for $\FL$.

Like the other results in the present paper, the claim P1 is only local:
it holds only for sufficiently small domain deformations.
A similar result is proved in WV with a precise functional setting and
a different (weaker) set of hypotheses.
The main idea, which we repeat here for reference, is to reformulate the
problem in terms of the area-preserving diffeomorphism
\begin{equation}\label{q:gL}
   \gL:D_\Lambdaz\to D_\Lambda:\xb\mapsto \gL\xb,
\end{equation}
which effects the vorticity rearrangement, that is, such that
$\wL = \omega_\Lambdaz \circ \gL^{-1}$. 
In terms of the pull-back 
\begin{equation}\label{q:g*def}
   \gL^*:\CFun(D_\Lambda)\to \CFun(D_\Lambdaz):f\mapsto f\circ\gL,
\end{equation}
this can be rewritten as
\begin{equation}\label{q:fund}
   \wL=(\gL^{-1})^*\omega_\Lambdaz.
\end{equation}
Note that, for fixed $\omega_\Lambdaz$ and $\wL$, $\gL$ is not defined
uniquely by \eqref{q:fund}: rearrangements along the lines of constant
vorticity clearly have no effect.
Correspondingly, the time derivative of $\gL$ is not necessarily the
divergence-free velocity field $\nabla^\bot \phi$ appearing in
\eqref{q:dw0dt}--\eqref{q:phidef}, but the equality
\begin{equation}\label{q:gphi}
   \ddt{}\gL \xb = \nabla^\bot [\phi(\gL \xb ; \Lambda)
				+ \varpi(\wL(\gL \xb)) ]
\end{equation} 
holds, where $\varpi$ is an arbitrary function of one variable.
This non-uniqueness, of no importance as far as $\psiL$ and $\wL$
are concerned, will play a crucial role when particle positions
are examined in~\S3.

The map $\gL$ satisfies a nonlinear partial differential equation
obtained as follows.
Since $\wL$ is a steady flow in $D_\Lambda$, we have using \eqref{q:Psi0},
\begin{equation}
   \wL = \lapl(\GL\circ\wL),
\end{equation}
so applying $\gL^*$, we find
\begin{equation}\label{q:glapg}
   \omega_\Lambdaz = \gL^*\lapl (\gL^{-1})^* (G_\Lambda\circ\omega_\Lambdaz).
\end{equation}
The associated boundary conditions are 
\begin{equation}
   \gL(\dy D_\Lambdaz) = \dy D_\Lambda
   \quad\textrm{and}\quad
   G_\Lambda = G_\Lambdaz 
\end{equation}
for the boundary value of $\omega_0$, 
the latter following from the fact that
$\psiL=\psi_\Lambdaz =0$ and $\wL=\omega_\Lambdaz$ on $\partial D_\Lambda$.

The partial differential equation \eqref{q:glapg}, with $\gL$ and $G_\Lambda$
as unknowns, is shown in WV to have a locally unique solution (modulo
translations along vorticity contours) using a contraction mapping argument.
This establishes P1 and provides a way of computing $\gL$ and $G_\Lambda$,
and hence $\wL$ and $\psiL$. 
Alternatively, P1 can be established using the stability assumption H1:
the associated characterisation of steady flows as energy extrema makes
it clear that the steady flow $\psiL$ is the (locally unique)
extremum in $D_\Lambda$ with vorticity distribution fixed by
$\omega_\Lambdaz$.

We now consider the infinitesimal version of \eqref{q:glapg}, that is,
we consider the change in $\gL$ corresponding to an infinitesimal
deformation of the domain.
This yields a different construction for $\gL$, based on integration
over $\Lambda$, and provides all the ingredients needed for the
computation of the first-order correction to $\psiL$ and of the
Lagrangian flow.

% ---------------------------------------------------------------------------

\subsection{Infinitesimal deformations}

In what follows, we will often make use of the fact that many important
quantities are linear in the boundary deformation rate
$\Lambdadot:={\rm d}\Lambda/{\rm d}\tau$.
We will regard these as resulting from the pairing between the vector
$\Lambdadot \in T_\Lambda \mathcal{L}$ and a differential one-form
belonging to a dual space.
For instance, the function $b$ appearing in the boundary condition
\eqref{q:bc0} is linear in $\Lambdadot$.
This makes it possible to define the one-form
\begin{equation} \label{q:nottoosure}
   \Psibc(\cdot;\Lambda):T_\Lambda\Lspace\to\CFun(\dy D_\Lambda)
\end{equation}
by 
\begin{equation} \label{q:bbeta}
   b = \Psibc\cdot\Lambdadot,
\end{equation}
where $\cdot$ denotes the pairing between vectors and one-forms.
Working with differential forms of this type gives a compact notation,
factoring out the dependence in $\Lambdadot$.
At the same time, it allows for a geometric interpretation of our results
as explained in the Discussion.

Let $\df$ be the exterior derivative in $\Lspace$.
Since
\begin{equation}
   \frac{\mathrm{d}}{\mathrm{d}\tau}{\gL} = \df \gL \cdot \Lambdadot, 
\end{equation}
and $\gL$ is area preserving,
we can define a function-valued one-form $\phiform$ by
\begin{equation}\label{q:dgL}
   \df\gL = \sgb\phiform\circ\gL.
\end{equation}
An equivalent statement is that
\begin{equation}\label{q:dgL2}
   \df\gL^* f  = [\phiform,f]
\end{equation}
for any $\Lambda$-independent function $f$.
The initial domain $\Lambda_0=\Lambdaz$ and initial flow flow $\psi_\Lambdaz$
having been fixed,
$\phiform$ depends only on $\Lambda$ and takes its value in the space of
functions in $D_\Lambda$; explicitly, 
\begin{equation}
   \phiform (\cdot;\Lambda) :T_\Lambda\Lspace\to\CFun(D_\Lambda).
\end{equation}
Alternatively, since $D_\Lambda \subset \RReal^2$, we can think of
$\phiform$ as a map from $\RReal^2$ to the cotangent bundle of
$\Lspace$, that is,
\begin{equation}
   \phiform : \RReal^2 \to T^*\Lspace.
\end{equation}

From \eqref{q:dgL}, $\phiform$ can be recognized as a connection one-form
encoding the change in $\gL$ that results from an infinitesimal change
in $\Lambda$ (making use of the identification between functions and
divergence-free vector fields on $D_\Lambda$; see the Discussion for a more
precise interpretation).
Introducing (formal) coordinates $\{\Lambda_m\}$ in (a neighbourhood
of $\Lambda=\Lambdaz$ in) $\Lspace$,  $\Phi$ takes the more explicit
form $\phiform=\phiform_m\,\df \Lambda_m$, where each $\phiform_m$
is a function in $D_\Lambda$ (sum over repeated indices is implied here
and henceforth). 

Taking the exterior derivative of \eqref{q:fund} gives
\begin{equation}\label{q:dfwL}
   \df\wL + [\phiform,\wL] = 0,
\end{equation}
after using the definition \eqref{q:dgL}.
In components,  this reads
\begin{equation}
   \ddy{\wL}{\Lambda_m} + [\phiform_m,\wL] = 0.
\end{equation}
Unlike $\wL$, the leading-order streamfunction $\psiL$ is not simply
rearranged as $\Lambda$ changes.
Applying $\df$ to $\psiL = \GL (\wL)$ and using \eqref{q:dfwL} show that
\begin{equation}
   \df\psiL + [\phiform,\psiL] = \df\GL\circ\wL.
\end{equation}
Here and elsewhere, in $\df\GL\circ\wL$ the exterior derivative is taken
with respect to the (parametric) dependence of $\GL$ on $\Lambda$,
so $\df(\GL\circ\wL) = \df\GL\circ\wL + (\GL'\circ\wL)\, \df\wL$.

From \eqref{q:dfwL} we derive a property of $\phiform$ which will be useful
later.
Taking $\df$ of \eqref{q:dfwL} and using the fact that $\df^2=0$,
leads, after a short computation detailed in Appendix~\ref{a:details}, to
\begin{equation}\label{q:integr}
\df\phiform + \sfrac12 [\phiform\wedge\phiform]
	= w\circ\wL 
\end{equation}
for some two-form $w\circ\wL:T^2_\Lambda\Lspace\to\CFun(D_\Lambda)$.
Here, the bracket $[\cdot\wedge\cdot]$ is defined in coordinates by
\begin{equation} \label{q:wbra}
[\alpha\wedge\gamma]=[\alpha_m,\gamma_n]\,\dLwdL{m}{n}
\end{equation}
for any two one-forms $\alpha=\alpha_m\,\df\Lambda_m$ and
$\gamma=\gamma_n\,\df\Lambda_n$ with values in $\CFun(D_\Lambda)$.
(Note that in contrast with $[f,f]=0$ for any function $f$,
$[\alpha\wedge\alpha]\ne0$ in general.)
Equation \eqref{q:integr} can be interpreted as an integrability condition,
an infinitesimal version of the statement that $\gL$ is a unique function
of $\Lambda$ modulo displacements along lines of constant $\wL$.  

We next obtain a dynamical equation for $\Phi$.
Taking the derivative $\df$ of \eqref{q:glapg} and after a little algebra,
we find this equation in the form
\begin{equation}\label{q:pde0}
   (\lapl-\FL')\,[\phiform,\psiL] - \lapl(\df\GL\circ\wL) = 0.
\end{equation}
The corresponding boundary conditions are
\begin{equation}\label{q:phibc}
   \phiform = \Psibc \qquad\textrm{on }\dy D_\Lambda,
\end{equation}
which follows from \eqref{q:bc0}, \eqref{q:bbeta} and \eqref{q:phidef},
taking into account that
\begin{equation}\label{q:chibc}
   \df\GL\circ\wL = 0 \qquad\textrm{on }\dy D_\Lambda,
\end{equation}
which follows from the fact $\psiL=0$ on $\dy D_\Lambda$.

Equation \eqref{q:pde0} is the infinitesimal version of \eqref{q:glapg}.
It can be solved for $[\phiform,\psiL]$ and $\df\GL$ as follows.
For any function $u \in \mathcal{C}(D_\Lambda)$, we define the projection
$\PpsiL u$ by
\begin{equation}
   \PpsiL\, u := u - \biggl\{ \oint_{\psiL=c} u \; \ds \biggr\} \bigg/
	\biggl\{ \oint_{\psiL=c} \ds \biggr\},
\end{equation}
where, as before, ${\rm d}s={\rm d}l/|\gb\psiL|$. 
It follows from this definition that if $u$ is constant on a contour
of constant $\psiL$, $\PpsiL u=0$. Note that since the contours of
$\wL$ and $\psiL$ coincide, a equivalent definition of $\PpsiL$ could
have been given in terms of integrals along vorticity contours $\wL=c$.
The regularity of $\PpsiL$ is guaranteed by H2.
Letting $\vfi=[\psiL,\phiform]-\df\GL\circ\wL$, and
using the facts that $\PpsiL[\psiL,\phiform]=[\psiL,\phiform]$
and $\PpsiL(\df\GL\circ\wL)=0$, we then have
\begin{equation}\label{q:phichivfi}
   [\psiL,\phiform] = \PpsiL\,\vfi
   \quad\textrm{and}\quad
   \df\GL\circ\wL = (1-\PpsiL)\,\vfi.
\end{equation}
Hence we can write \eqref{q:pde0} as
\begin{equation}\label{q:vfi}
   (\Delta - \FL'\,\PpsiL)\,\vfi = 0,
\end{equation}
which is a linear pseudodifferential equation involving $\vfi$ only.
Following \eqref{q:phibc} and \eqref{q:chibc}, the boundary conditions
for $\vfi$ are
\begin{equation}\label{q:vfibc}
   \vfi = [\psiL,\Psibc]
   \qquad\textrm{on }\dy D_\Lambda.
\end{equation}
It is shown in Appendix~\ref{a:lin} that \eqref{q:vfi}--\eqref{q:vfibc} can be
solved uniquely for $\vfi$.
Using \eqref{q:phichivfi}, we recover $[\psiL,\phiform]$ and $\df\GL$.
From $[\psiL,\phiform]$, $\phiform$ can be inferred up to an arbitrary
function-valued one-form depending on $\xb$ through $\psiL$ or,
equivalently, through $\wL$.
Thus, there is an equivalence class of forms $\phiform$ satisfying
\eqref{q:pde0} which is associated with the gauge transformation
\begin{equation}\label{q:gauge}
   \Phi \mapsto \Phi + \Pi \circ \wL,
\end{equation}
where $\Pi \circ \wL$ is any function-valued one-form depending on $\xb$
through $\wL$.
This non-uniqueness simply reflects at the infinitesimal level the
non-uniqueness of $\gL$.
A particular, uniquely-defined representative of the equivalence class
of $\Phi$ could be taken to be $\PpsiL \Phi$, that is, the one with
vanishing average along streamlines $\psiL=\textrm{const}$.
This is an arbitrary choice, however, and we will see in the next
section that another choice imposes itself naturally. 

As shown in WV, once we solve the linear problem \eqref{q:pde0},
the solution of the nonlinear problem \eqref{q:glapg} for $\gL$ follows,
at least in a neighbourhood of $\Lambda=\Lambdaz$ and subject to sufficient
smoothness of the flow and the domain.
For a fixed sequence of boundary deformation, that is, for a fixed path
$\Lambda(\tau)\subset\Lspace$, one can in principle solve \eqref{q:pde0}
and find $\phiform$ for each $\Lambda$ as long as the flows
$\psiL$ encountered along the path satisfy H1--H2.
Note that \eqref{q:pde0} is consistent with the integrability condition
\eqref{q:integr} which thus remains satisfied along the path (see
Appendix~\ref{a:details}).
This confirms our main conclusion, namely that $\gL$ is independent of
the path, up to translation along contours of $\wL$. 

% ===========================================================================

\section{Eulerian Flow: First-order Correction}\label{s:eu1}

We now turn to the derivation of the first correction $\psio$
to the leading-order flow $\psiL$.
This derivation  is necessary, in particular, to determine the trajectories of
fluid particles over the $\Oh(\eps^{-1})$ time scales of interest. 
Remarkably, $\psio$ can be derived from the knowledge of $\phiform$ alone.
In the process, the gauge of $\phiform$ is fixed in what we argue is a natural
manner.

The derivation starts by noting that \eqref{q:gphi} and \eqref{q:dgL}
imply that $\phi$ and $\Phi \cdot  \Lambdadot$ differ by a function
of $\wL$ only.
Thus we can write
\begin{equation} \label{q:phiphiform}
   \phi = \phinat \cdot \Lambdadot,
\end{equation}
where
\begin{equation} \label{q:phinatphiform}
 \phinat = \PpsiL \phiform + \Pinat \circ \wL.
\end{equation}
Here $\phinat$ is a specific member of the class of equivalent one-forms
$\Phi$: it corresponds to the unique choice of the gauge $\Pi = \Pinat$
in \eqref{q:gauge} that ensures that \eqref{q:phiphiform} holds.
In this sense $\phinat$ can be seen as a natural choice of connection form.
The  computation which follows shows how it can be obtained from
$\PpsiL\phiform$.

Since $\psio$, like $\phi$, is linear in $\Lambdadot$, we can write
\begin{equation}
   \psio=\Psio \cdot \dot \Lambda,
\end{equation}
where, $\Psio$, like $\phiform$, is a function-valued form;
their relationship follows from (\ref{q:phidef}) as
\begin{equation} \label{q:PsiPi}
   (1-\GL'\Delta) \Psio = \PpsiL \phiform + \Pinat \circ \wL.
\end{equation}
Both $\Psio$ and $\Pinat$ can be deduced from (\ref{q:PsiPi}).
To show this, we make use of the constraint imposed by the material
conservation of the total vorticity $\omega$.
This implies that the function
\begin{equation}
   \Area(\Omega;\omega)
	:= \int_{\textrm{int}\{\omega=\Omega\}} \d2x
	= \textrm{area bounded by }\{\omega=\Omega\},
\end{equation}
where int$\{\omega=\Omega\}$ denotes the interior of the curve $\omega=\Omega$,
is an exact invariant of the dynamics.
Since both the total and leading-order vorticities $\omega$ and $\wL$
are rearrangements of the initial vorticity $\omega_{\Lz}$,
\begin{equation}
   \Area(\Omega;\omega) =  \Area(\Omega;\wL) = \Area(\Omega;\omega_{\Lz}).
\end{equation}
Expanding $\omega=\wL + \eps \omegao + \cdots$, we then have
\begin{equation}
   \Area(\Omega;\wL+\eps\omegao) - \Area(\Omega;\wL) = \Oh(\eps^2),
\end{equation}
or, after some manipulations,
\begin{equation}\label{q:w1gauge}
   \oint_{\wL=\Omega} \omegao \ds = 0.
\end{equation}
The last equation, which can be rephrased as $(1-\PpsiL)\,\omegao = 0$ or
\begin{equation} \label{q:PPsi1}
   (1-\PpsiL) \Delta \Psio = 0,
\end{equation}
provides the solvability condition for \eqref{q:PsiPi}.
Indeed, applying $(1-\PpsiL)$ to \eqref{q:PsiPi} gives
\begin{equation}\label{q:PiPsio}
   \Pinat \circ \wL = (1-\PpsiL) \Psio.
\end{equation}
This reduces \eqref{q:PsiPi} to 
\begin{equation}\label{q:psi1}
   (\Delta - \FL' \PpsiL) \Psio = - \FL' \PpsiL \phiform.
\end{equation}
The associated  boundary condition follows from  \eqref{q:psibc} as
\begin{equation} \label{q:Psiobc}
  \Psio = \Psibc \quad \textrm{on }\> \partial D_\Lambda.
\end{equation}

Equation~\eqref{q:psi1} is well posed, with a right-hand side that is uniquely
defined in spite of the gauge freedom in $\phiform$.
The operator on the left-hand side is the same as that in (\ref{q:vfi}) and
hence its invertibility can be established using the same arguments,
detailed in Appendix~\ref{a:lin}.
Once $\Psio$ is determined from \eqref{q:psi1}, $\Pinat$ follows from
(\ref{q:PiPsio}), and the natural connection $\phinat$ is obtained. 
Note that \eqref{q:phinatphiform} and \eqref{q:PiPsio} imply that it satisfies
\begin{equation}\label{q:psi1mean}
   (1-\PpsiL) \Psio = (1-\PpsiL) \phinat
   \qquad\Leftrightarrow\qquad
   \oint_{\psiL=c} \Psio \ds = \oint_{\psiL=c} \phinat \ds.
\end{equation}
This relation turns out to be crucial for the computation of fluid particle
trajectories in the next section.
  
% ===========================================================================

\section{Lagrangian Flow: Geometric Angle}\label{s:lag}

In this section we study the evolution of fluid (or tracer) particles in
our flow over a timescale $\tau = \Oh(1)$.

% ---------------------------------------------------------------------------

\subsection{Hamiltonian Formulation}

Up to this point, our description of the Eulerian dynamics has been
(mostly) coordinate-independent.
But in order to study particle positions, we need to introduce explicit
coordinates $(x,y)$ in $D_\Lambda$;
$(x,y)$ is chosen to coincide with the fixed coordinates in the ambient
space $\RReal^2$ through which $D_{\Lambda(\eps t)}$ moves.

The evolution of a particle with position $(x(t),y(t))$ moving with the fluid
is governed by the Hamiltonian system
\begin{equation}\label{q:dxdt}
   \ddt{x} = -\ddy{\psi}{y}
   \quad\textrm{and}\quad
   \ddt{y} = \ddy{\psi}{x},
\end{equation}
with the streamfunction $\psi$ acting as the Hamiltonian.
Our aim here is to obtain an estimate of $(x(t),y(t))$ with an error of
$\Oh(\eps)$ for $\tau = \Oh(1)$, so in the rest of this section we put
\begin{equation}\label{q:h0h1}
   H(x,y,t) = \psi_{\Lambda(\eps t)}(x,y)
		+ \eps\Psio(x,y;\Lambda(\eps t))\cdot\Lambdadot
\end{equation}
in place of $\psi$ in \eqref{q:dxdt}, keeping in mind the validity of this
approximation. 

For $\eps=0$ and hence $\Lambda$ constant, the Hamiltonian \eqref{q:h0h1}
is integrable.
For $\eps \not=0$, two types of perturbations make it non-integrable:
the slow time dependence of $\psiL$ introduced by the time dependence
of $\Lambda$, and the $\Oh(\eps)$ change introduced by the addition of
$\psio=\Psio \cdot \Lambdadot$.
We examine the combined effect of these two perturbations following
closely the approach of \cite{berry:85}.

Since the leading-order Hamiltonian $\psiL$ is integrable for fixed $\Lambda$,
we first change to action--angle variables
[cf.\ \cite[pp.~297ff]{arnold:mmcm}].
At each $(x,y)$, we define the action $I$ by
\begin{equation}\label{q:Idef}
   I(x,y) = \frac{1}{2\pi} \int_{\textrm{int}\{\psiL=\psiL(x,y)\}} \dxdy
	=: \frac{1}{2\pi} A(\psiL).
\end{equation}
The angle $\theta$ is defined as the variable conjugate to $I$, $[I,\theta]=1$.
It is $2\pi$-periodic since the contours of $\psiL$ are closed
and it is related to the variable $s$ used earlier by
\begin{equation} \label{q:thetas}
 2\pi\,\mathrm{d}s = A'(\psiL)\, \mathrm{d}\theta.
\end{equation}
The canonical transformation $(x,y)\mapsto(I,\theta)$ is obtained
by a generating function $S(I,y;\Lambda)$, with
\begin{equation}
   x = \ddy{S}{y}
   \quad\textrm{and}\quad
   \theta = \ddy{S}{I}.
\end{equation}
Solving these implicit equations, we can write
\begin{equation}
   x = X(I,\theta;\Lambda)
   \quad\textrm{and}\quad
   y = Y(I,\theta;\Lambda).
\end{equation}
With these and \eqref{q:Idef}, we define
\begin{equation}\label{q:hat}\begin{aligned}
   &\psizb(I;\Lambda) = \psiz(X(I,\theta;\Lambda),Y(I,\theta;\Lambda);\Lambda)
	= A^{-1}(2\pi I;\Lambda),\\
   &\Psioh(I,\theta;\Lambda)
	= \Psi^{(1)}(X(I,\theta;\Lambda),Y(I,\theta;\Lambda);\Lambda).
\end{aligned}\end{equation}
Here and in the rest of this section, we denote by a h{\^a}t quantities
considered as functions of $(I,\theta)$.

So far we considered a fixed value of $\Lambda$.
Now let $\Lambda$ evolve slowly in time, $\Lambda=\Lambda(\eps t)$.
The equations of motion in $(I,\theta)$ variables are
\begin{equation}\label{q:dIthdt}
   \ddt{I} = -\ddy{\hat H}{\theta}
   \quad\textrm{and}\quad
   \ddt{\theta} = \ddy{\hat H}{I},
\end{equation}
where the new Hamiltonian $\hat H(I,\theta;\Lambda)$ is related to
$H(x,y;\Lambda)$ by
\begin{equation}\label{q:Hhat}
   \hat H(I,\theta;\Lambda)
	= H(X(I,\theta;\Lambda),Y(I,\theta;\Lambda);\Lambda) + \ddy{S}{t}.
\end{equation}
We note an abuse of notation here: properly speaking $H=H(x,y,t)$,
but since the $t$-dependence only enters through $\Lambda(\eps t)$ and its
derivative, we have written $H=H(x,y;\Lambda(\eps t))$.
Differentiating the definition
\begin{equation}
   \hat S(I,\theta;\Lambda) = S(I,Y(I,\theta;\Lambda);\Lambda)
\end{equation}
with respect to $t$ at fixed $(I,\theta)$ gives
\begin{equation}
   \df \hat S \cdot\ddt\Lambda =  \ddy{S}{t} + \ddy{S}{y}\ddy{Y}{t}
		= \ddy{S}{t} + X \df Y \cdot \ddt\Lambda
\end{equation}
Upon substituting $\dy S/\dy t$ into \eqref{q:Hhat}, we obtain that
\begin{equation}\begin{aligned}
   \hat H(I,\theta;\Lambda) = \psizb(I)
	&+ \eps\Psioh(I,\theta;\Lambda)\cdot\Lambdadot\\
	&+ \eps \bigl\{ \df{\hat S}(I,\theta;\Lambda)
	   - X(I,\theta;\Lambda) \df{Y}(I,\theta;\Lambda) \bigr\}
	\cdot\Lambdadot.
\end{aligned}\end{equation}
Since particles are attached to contours of vorticity $\omega=\textrm{const}$,
which only deviate by $\Oh(\eps)$ from the corresponding contours of $\wL$,
the action can only vary by $\Oh(\eps)$ over timescales $\tau\sim\Oh(1)$.
This is also evident from direct computation:
since $\hat H$ is independent of $\theta$ at leading order and is periodic
in $\theta$ at the next order,
\begin{equation}
   \ddt{I} = -\eps\, \ddy{}{\theta}\bigl\{ \Psioh + \df{\hat S}
	   - X \df{Y} \bigr\} \cdot \Lambdadot,
\end{equation}
and the principle of averaging [cf.\ \cite[\S~52]{arnold:mmcm}]
implies that $I$ changes only by $\Oh(\eps)$ for $\tau = \Oh(1)$.

The behaviour of the angle variable is more interesting.
From \eqref{q:dIthdt} we have
\begin{equation}\label{q:dthdt}
   \ddt{\theta} = \ddy{\psizb}{I}
	+ \eps\, \ddy{}{I} \bigl\{ \Psioh + \df{\hat S} - X \df{Y} \bigr\} 
	\cdot\Lambdadot.
\end{equation}
The change in the angle $\Delta\theta:=\theta(\tau)-\theta(0)$ can then be
expressed as $\Delta\theta = \dynp + \geop$.
The dynamic phase $\dynp$ simply arises from the instantaneous
frequency of the particle, which is the first term in \eqref{q:dthdt} above,
\begin{equation}
   \dynp = \frac{1}{\eps} \ddy{}{I} \int_0^\tau \hat\psi_{\Lambda(\tau')}(I)
	 \;{\rm d}\tau'.
\end{equation}
The other terms make up the geometric angle $\geop$.

% ---------------------------------------------------------------------------

\subsection{Geometric angle $\geop$}

In this subsection we show that, as in the
finite-dimensional cases of \cite{hannay:85} and \cite{berry:85},
the angle $\geop$ can be understood in geometric terms as the
(an)holonomy of a connection as a closed path is traversed
in a parameter space.

From \eqref{q:dthdt}, the geometric angle $\geop$ can be written as
\begin{equation}\label{q:geop}\begin{aligned}
   \geop = \int_0^\tau \ddy{}{I}\bigl\{ \Psioh &(I,\theta;\Lambda(\tau'))
		+ \df{S}(I,\theta;\Lambda(\tau'))\\
	&- X(I,\theta;\Lambda(\tau'))\,\df{Y}(I,\theta;\Lambda(\tau'))\bigr\}
		\cdot\frac{{\rm d}\Lambda}{{\rm d}\tau'} \;{\rm d}\tau'
\end{aligned}\end{equation}
This form suggests that $\geop$ depends only on the path traversed
in $\Lspace$ and not on its time parametrisation.
The terms inside the braces do depend on $I$ and $\theta$, but as shown
earlier, the total variation of the action $I$ is of $\Oh(\eps)$ over
the timescale of interest.
The dependence on the periodic variable $\theta$ can be removed by
averaging.
For any function $f$ periodic in $\theta$, let
\begin{equation}
   \langle f\rangle := \frac{1}{2\pi} \int_0^{2\pi} f(\theta) \;{\rm d}\theta.
\end{equation}
We note that by \eqref{q:thetas} this is essentially equivalent to the
projection $1-\PpsiL$.
Applying $\langle\cdot\rangle$ to \eqref{q:geop} and replacing $I(t)$ by
$I(0)$, we find
\begin{equation}\label{q:geop1}
   \geop = \ddy{}{I} \int_{C_\Lambda} \bigl\langle \Psioh
	+ \df S - X\,\df Y \bigr\rangle + \Oh(\eps),
\end{equation}
where $C_\Lambda$ is the path traversed in $\Lspace$ and where
the integrand depends only on $I$ and $\Lambda$.

\begin{figure}
\begin{center}
\epsfig{file=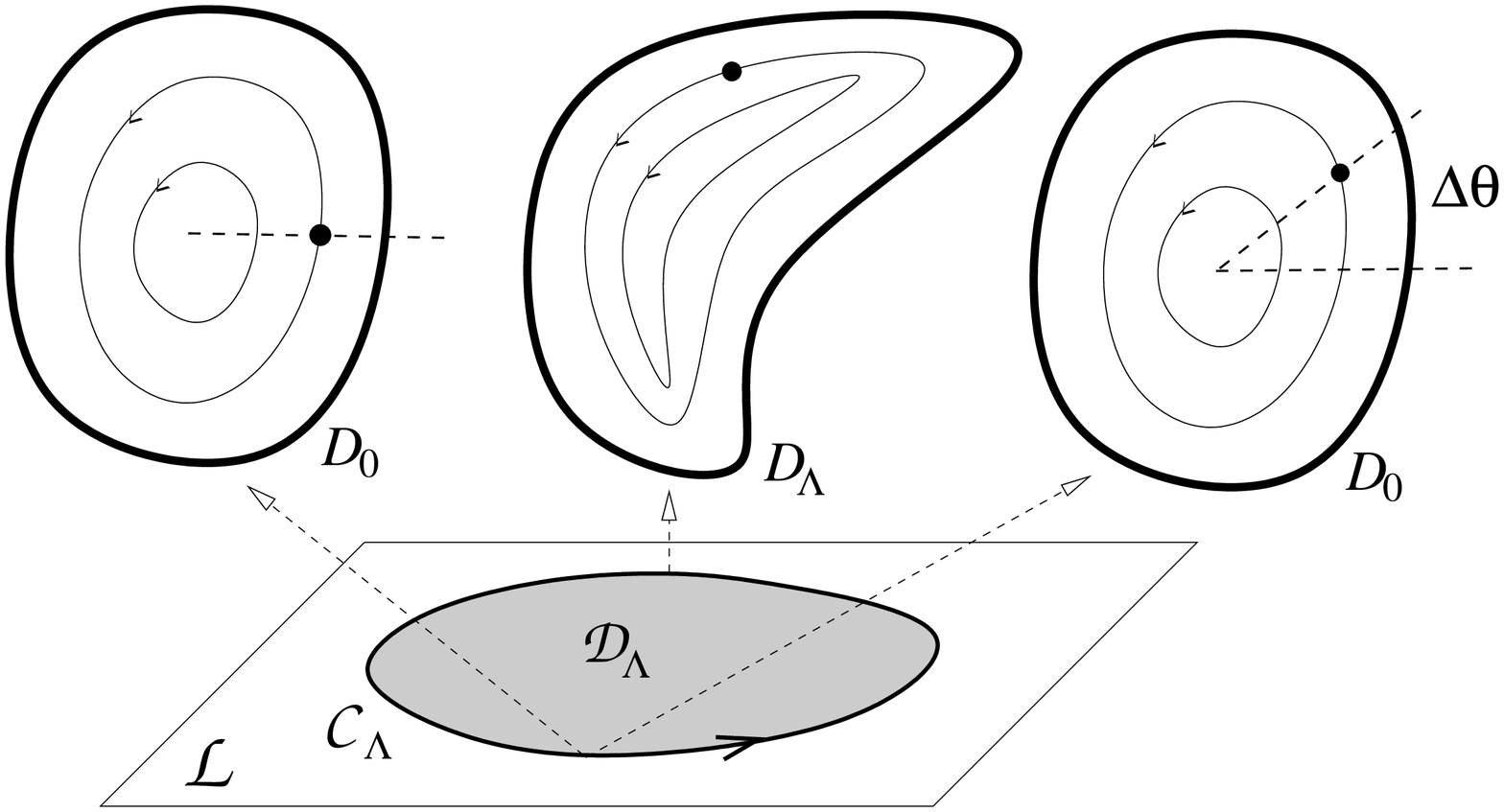,height=5cm}
\caption{%
Angle change for a cyclic domain deformation. As $\Lambda(\tau)$
describes the closed loop $\mathcal{C}_\Lambda \subset \mathcal{L}$,
with interior $\mathcal{D}_\Lambda$, the fluid domain $D_\Lambda$ is
deformed and returns to its original shape.
Fluid particles remain on vorticity contours which approximately coincide
with streamlines. The position of the particles along vorticity contours is
defined by the angle-like variable $\theta$ whose total change $\Delta\theta$
includes the geometric contribution $\Delta\theta_{\mathrm{geo}}$ which
depends only on the geometrical properties of $\mathcal{D}_\Lambda$.
}\label{fig:loop}
\end{center}
\end{figure}

Because of the arbitrariness in the angle coordinates (depending on our
choice of $\theta=0$ for each $\Lambda$), the geometric angle is only
unambiguously defined when the path $C_\Lambda$ is closed, that is, when
$\Lambda(\tau)=\Lambda(0)=0$.
Following \cite{hannay:85} and \cite{berry:85}, we consider this scenario,
which is illustrated in Figure \ref{fig:loop}.
Since $\df S$ is exact, it vanishes when integrated around $C_\Lambda$.
Using Stokes' theorem in $\Lspace$, the remaining terms in \eqref{q:geop1}
can be written as
\begin{equation}\label{q:geop2}
   \geop = \ddy{}{I} \int_{\mathcal{S}_\Lambda} \bigl\langle
	\df\Psioh - \df X\wedge\df Y \bigr\rangle
\end{equation}
where $\mathcal{S}_\Lambda$ is a two-dimensional surface bounded
by $C_\Lambda$.
The second term is identical to that obtained by Berry (1985; eq.~(18))
for general Hamiltonian systems
depending slowly on time;
the first term results from the $\Oh(\eps)$ change to the Hamiltonian
induced by the boundary deformation.

\nocite{berry:85}

Now $\langle\Psioh\rangle = \langle\phinath\rangle$ by \eqref{q:psi1mean},
so the first term in the integral can be written as
\begin{equation}
   \langle\df\Psioh\rangle = \langle\df\phinath\rangle.
\end{equation}
As for the second term, we use the fact that $\gL$ is a canonical
transformation (since it is area-preserving) to write
\begin{equation} \label{q:X}
   \Xb(I,\theta;\Lambda) = \gL\Xb(I,\theta;\Lambdaz),
\end{equation}
thus defining the transformation to action--angle coordinates for all
values of $\Lambda$ in terms of the transformation at $\Lambda=\Lambdaz$.

It follows from this and \eqref{q:dgL} that
\begin{equation}
   \df\Xb(I,\theta;\Lambda) = \sgb\phiform\big|_{\Xb(I,\theta;\Lambda)}.
\end{equation}
We then have
\begin{equation}\begin{aligned}
   \df X\wedge\df Y
	&= \ddy{X}{\Lambda_m}\ddy{Y}{\Lambda_n} \;\dLwdL{m}{n} \hspace{-20pt}
	&&= -\ddy{\phiform_m}{Y}\ddy{\phiform_n}{X} \;\dLwdL{m}{n}\\
	&= \sfrac12 [\phiform_m,\phiform_n] \;\dLwdL{m}{n}
	&&= \sfrac12 [\phiform\wedge\phiform]
	= \sfrac12 [\phinat\wedge\phinat],
\end{aligned}\end{equation}
where the last equality follows from the fact that the bracket
is independent of the gauge choice for $\phiform$.
Furthermore, the fact that the transformation to action--angle
variables is canonical implies that
\begin{equation}
   \sfrac12 [\phinat\wedge\phinat] = \sfrac12 [\phinath\wedge\phinath],
\end{equation}
where the second bracket is in terms of $(I,\theta)$. 
We can therefore write \eqref{q:geop2} as
\begin{equation}\label{q:geop3}
   \geop = \frac{{\rm d}}{{\rm d}I} \int_{\mathcal{S}_\Lambda} \bigl\langle
	\df\phinath - \sfrac12[\phinath\wedge\phinath] \bigr\rangle.
\end{equation}
A last step, detailed in Appendix~\ref{a:details}, shows that 
\begin{equation} \label{q:kappa}
   \kappa := \df\phinath - \sfrac12[\phinath\wedge\phinath]
	  = \df\phinat + \sfrac12[\phinat\wedge\phinat].
\end{equation}
This function-valued two-form can be recognised as the curvature of
the connection $\phinat$, and according to \eqref{q:integr},
it depends on space through $\wL$ or, equivalently, through $I$ only.
This property stems from  the constraints of area and vorticity preservation
imposed on the particle motion along vorticity contours.
The average in \eqref{q:geop3} is therefore superfluous,
and we obtain the result:

\medskip\noindent{\bf P2.}
{\sl The geometric angle of a particle caused by the slow deformation
of the boundary is given by
\begin{equation}\label{q:geop4}
   \geop = \frac{{\rm d}}{{\rm d}I} \int_{\mathcal{S}_\Lambda} \kappa,
\end{equation}
where $\kappa$, given in \eqref{q:kappa}, is the curvature of $\phinat$ and
depends only on the action $I$ of the particle and on the domain shape
parameterised by $\Lambda$.} 

% ===========================================================================

\section{Small boundary deformation}\label{s:small}

As a concrete illustration of the developments so far, we consider the
case where the total boundary deformation is small.
In general, computing $\gL$ in \eqref{q:gL} requires either solving
\eqref{q:glapg} or integrating the differential equation \eqref{q:dw0dt}
with boundary condition $\phi=b$ [cf.\ \eqref{q:phibc}],
which would have to be done numerically.
Analytic progress is possible, however, if one considers boundary
deformations that are sufficiently small for a perturbative approach
to be applicable.
In this section we develop such an approach systematically using Lie series
[cf.\ \cite{lichtenberg-lieberman:rcd}]. 

Let $\delta$ be a formal small parameter characterising the smallness
of the boundary deformation.
The function $B(\xb;\Lambda)$ defining $D_{\Lambda}$ can then be expanded as
\begin{equation}
   B(\xb;\Lambda) = B_{0}(\xb) + \delta B_{1}(\xb;\Lambda)
	+ \delta^{2}  B_{2}(\xb;\Lambda) + \cdots,
\end{equation}
where $B_{0}$ is independent of $\Lambda$ (recall that $\Lambda_0=\Lambdaz$),
$B_{1}$ is linear in $\Lambda$,
$B_{2}$ quadratic, etc. 
Since $\gL$ is close to the identity and area-preserving, it may be
regarded as the flow at `time' $\delta$ of an associated $\delta$-dependent
divergence-free vector field given by $\nabla^\bot \rho$ for some function
$\rho(\xb,\delta)$.
(Here $\rho$ is to be considered to live in $\CFun(ND\times I)$, where
$ND\subset\RReal^2$ is large enough to contain all relevant $D_\Lambda$
and $I\subset\RReal$.)
Correspondingly, the pull-back of $\gL$ defined in \eqref{q:g*def} satisfies
\begin{equation}\label{q:liederiv}
   \frac{\mathrm{d}\gL^{*}}{\mathrm{d}{\delta}} = \gL^{*} [\rho, \cdot].
\end{equation}
Expanding  $\rho$ in powers of $\delta$ as
\begin{equation}
   \rho = \rho_{1} + \delta \rho_{2} + \cdots
\end{equation}
and introducing into (\ref{q:liederiv}) lead to the expansions
\begin{align}
   \gL^{*} &= 1 + \delta [\rho_{1}, \cdot]
	+ \frac{\delta^{2}}{2} \bigl([\rho_2,\cdot]
		+ [\rho_1,[\rho_1,\cdot]] \bigr) + \cdots, \label{q:g*}\\
   (\gL^{-1})^{*} &= 1 - \delta [\rho_{1}, \cdot]
	- \frac{\delta^{2}}{2} \bigl(  [\rho_2,\cdot]
		- [\rho_1,[\rho_1,\cdot]] \bigr) + \cdots. \label{q:g-*}
\end{align}
Computing $\gL^*f$ for an arbitrary $\Lambda$-independent $f$ using
\eqref{q:g*}, taking the exterior derivative and identifying with
\eqref{q:dgL2} lead to
\begin{equation}\label{q:Phirho}
   \phiform = \delta \df \rho_1+ \frac{\delta^{2}}{2}
	\bigl( \df \rho_2 + [\df\rho_1,\rho_1] \bigr) + \cdots,
\end{equation}
up to an arbitrary function of $\wL$, after using the Jacobi identity. 

Introducing \eqref{q:g*}--\eqref{q:g-*} into (\ref{q:glapg}) leads to a
sequence of partial differential equations for the coefficients of $\rho$. 
The first two read
\begin{align}
   (\Delta - F_{\Lz}') [\rho_{1},\psi_{\Lz}] - \Delta \chi_{1}^{} &= 0,
	\label{q:varphi1} \\
   (\Delta - F_{\Lz}') [\rho_{2},\psi_{\Lz}] - 2 \Delta \chi_{2}^{}
	&= -2 [ \rho_1,\Delta[ \rho_1,\psi_{\Lz}]]
	+ [\rho_1,[\rho_1,\omega_{\Lz}]]   \label{q:varphi2} \\
	&\qquad+ \Delta [\rho_1,[\rho_1,\psi_{\Lz}]]
	+ 2 [\rho_1,\Delta \chi_1^{}]-2 \Delta [\rho_1,\chi_1^{}], \notag
\end{align}
where $F_{\Lz}'$ is shorthand for $F_{\Lz}'\circ\psi_{\Lz}$, and we have
introduced the expansion
\begin{equation}
  \GL\circ\omega_{\Lz} = \psi_{\Lz} + \delta \chi_{1}^{}
	+ \delta^{2} \chi_{2}^{} \,, \cdots,
\end{equation}
with $\chi_i^{}$, $i=1,2\cdots$, depending on $\xb$ through
$\omega_{\Lz}(\xb)$.
These equations are supplemented by the boundary conditions
\begin{align}
   [\rho_1,B_0] &= - B_1,  \label{q:bc1} \\
   [\rho_2,B_0] &= - 2 B_2 - 2 [\rho_1,B_1] -
	[\rho_1,[\rho_1,B_0]], \label{q:bc2} 
\end{align}
to be applied on the curve $B_0(\xb)=B(\xb;\Lambdaz)=0$.
The formulation is then relatively simple, with all the equations to be
solved in the original domain $D_{\Lz}$.
The functions $\chi_{n}^{}$, $n=1,2,\cdots$ are found from solvability
conditions.
These can be made explicit using the same method as in the treatment
of \eqref{q:pde0}.
For instance, using the projection operator $\Ppsiz$ associated with
lines of constant $\psi_{\Lz}$, \eqref{q:varphi1} can be rewritten as
\begin{equation}\label{q:vfi1}
   (\lapl-F_\Lz'\Ppsiz)\,\vfi = 0
\end{equation}
where $\vfi:=[\rho_1,\psi_{\Lz}]-\chi_1^{}$, implying that
\begin{equation}
   [\rho_1, \psi_{\Lz}] = \Ppsiz \varphi
   \quad\textrm{and}\quad
   \chi_{1}^{} = (1-\Ppsiz) \varphi.
\end{equation}
Once the $\rho_{n}$, $n=1,2,\cdots$, are computed, the leading-order
vorticity and streamfunction follow readily from
\begin{align}
   \wL &= \omega_\Lz - \delta [\rho_1,\omega_\Lz] - \frac{\delta^2}{2}
   \bigl([\rho_2,\omega_\Lz] -
	[\rho_1,[\rho_1,\omega_\Lz]]\bigr)+\cdots, \\
   \psiL &= \psi_\Lz - \delta \bigl([\rho_1,\psi_\Lz] -\chi_1^{}\bigr)\notag\\
	&\qquad- \frac{\delta^2}{2} \bigl([\rho_2,\psi_\Lz] -
	  [\rho_1,[\rho_1,\psi_\Lz]]- 2 \chi_2^{}\bigr )+\cdots.
\end{align}

To find the first-order correction to the Eulerian flow, \eqref{q:psi1}
needs to be solved by expansion in powers of $\delta$.
This is conveniently done by pulling back this equation to the original
domain $D_{\Lz}$.
To do this, we define the pull-backs (denoted by overbars) and their
expansions as
\begin{align}
   \Psib1 &:= \gL^{*} \Psio
	& &\hspace{-25pt}= \delta \Psib1_1 + \delta^{2} \Psib1_2 + \cdots, \\ 
   \bar\phiform &:= \gL^{*} \phiform
	& &\hspace{-25pt}= \delta \bar{\phiform}_1
		+ \delta^{2} \bar\phiform_2 + \cdots \notag\\
	& & &\hspace{-25pt}= \delta \df \rho_{1} + \frac{\delta^{2}}{2}
		\bigl( \df \rho_{2} - [\df \rho_{1},\rho_{1}] \bigr)
	+ \cdots, \hbox to50pt{}
\end{align}
where the last equality follows from \eqref{q:g*} and \eqref{q:Phirho}.
Introducing these pull-backs into \eqref{q:psi1} and noting that
\begin{align*}
   (\gL^{-1})^* F'_\Lambda(\psiL)
	&= (\gL^{-1})^* [G'_\Lambda(\wL)]^{-1}
	= [G'_\Lambda(\omega_{\Lz})]^{-1} \\
        &= F'_{\Lz}(\psi_\Lz) - \delta [F'_{\Lz}(\psi_\Lz)]^2
		\frac{\nabla \chi_1^{}}{\nabla \omega_{\Lz}} + \cdots
\end{align*}
leads to
\begin{align}
   (\Delta - F'_{\Lz} \Ppsiz) \Psib1_1
	&=  - F'_{\Lz} \Ppsiz \bar\phiform_1, \label{q:noname1} \\
   (\Delta - F'_{\Lz} \Ppsiz) \Psib1_2
	&= - F'_{\Lz} \Ppsiz \bar\phiform_2
	+ \Delta [ \rho_{1},\Psib1_1] - 
	[ \rho_{1},\Delta \Psib1_1] \nonumber \\
	&\qquad-  F'_{\Lz} \frac{\nabla \chi_1^{}}{\nabla \omega_{\Lz}}
		\Delta \Psib1_1. \label{q:noname2}
\end{align}
These equations, involving the same invertible operator as
\eqref{q:vfi1}, can be solved to find $\Psib1$,
with $\Psio$ deduced after application of $(\gL^{-1})^{*}$.
The natural gauge $\phinat$ of $\phiform$ then follows from
\eqref{q:phinatphiform} and \eqref{q:PiPsio}.
Alternatively, one can first compute  the pull-back $\phinatb$,
which is obtained from the relations
\begin{equation} \label{q:phibs}
   \phinatb = \Ppsiz \bar \phiform + \Pi^\star \circ \omega_0
   \quad\textrm{and}\quad
   \Pi^\star \circ \omega_0 = (1-\Ppsiz) \Psiob
\end{equation}
inferred from \eqref{q:phinatphiform} and \eqref{q:PiPsio},
and then deduce $\phinat$ by pushing forward with $\gL^*$. 

To compute the curvature $\kappa$ and the geometric angle, there is in fact
no need to push forward $\Psib1$ and $\phinatb$: indeed, from
\eqref{q:hat} and \eqref{q:X}, we see that $\Psioh$ and $\Psib1$
are related by the $\Lambda$-independent transformation
\begin{equation}
   \Psioh(I,\theta;\Lambda)=\Psib1(\Xb(I,\theta;\Lambdaz);\Lambda)
\end{equation}
defining the action--angle variables in the original domain $D_\Lambdaz$. 
Since $\phinath$ and $\phinatb$ obey an analogous relation,
they are essentially equivalent: in particular, $\df \phinath = \df \phinatb$
and $[\phinath \wedge \phinath] = [\phinatb \wedge \phinatb]$.
The curvature $\kappa$ in \eqref{q:kappa} can therefore be computed directly
from $\phinatb$ in a straightforward manner as
\begin{equation}\label{q:kapphinat}
   \kappa = \df\phinatb - \sfrac12 [\phinatb\wedge\phinatb].
\end{equation}

Note that, in principle, the first two terms in the expansion of $\phinatb$
or $\phinat$ need to be computed in order to obtain a leading-order
approximation to the geometric angle.
This is because $\phinatb_{1}$ is independent of $\Lambda$,
$\df\phinatb_{1}=0$ and hence $\kappa = \Oh(\delta^{2})$.
The computation can however be shortened by observing that the average
of $\phinatb_{2}$ along streamlines, that is, $(1-\Ppsiz)\phinatb_{2}$,
is the only $\Oh(\delta^2)$ quantity genuinely needed if the averaged
form \eqref{q:geop3} of $\kappa$ is used.
In turn, $(1-\Ppsiz)\phinatb_{2}$ can be approximated by $(1-\Ppsiz)\Psiob_2$,
as the pull-back of \eqref{q:psi1} indicates.
The latter quantity satisfies a  relatively simple equation, obtained by
applying $(1-\Ppsiz)$ to \eqref{q:noname2} to find
\begin{equation} \label{q:avpsiob}
   (1 - \Ppsiz) \Delta \Psib1_2 = (1-\Ppsiz) \bigl\{ 
	\Delta [ \rho_{1},\Psib1_1] - 
	[ \rho_{1},\Delta \Psib1_1] \bigr\},
\end{equation}
after using $(1-\Ppsiz) \Delta \Psib1_1=0$ which follows from
\eqref{q:PPsi1} at leading-order in $\delta$. 

% ===========================================================================

\section{Nearly axisymmetric flow}\label{s:axi}

We now consider a simple example where the computations of $\gL$ and other
relevant quantities can be carried out explicitly to $O(\delta^{2})$.
We assume that for $\Lambda=\Lz$, the fluid domain is the disc
$(r,\sig)\in[0,1] \times [0,2 \pi]$.
The deformed domain is defined by 
\begin{equation} \label{q:disc}
   r = 1 + \delta \sum_{m} \Lambda_m \ex^{\i m \sig}
	- \frac{\delta^2}{2} \sum_m |\Lambda_m|^2 + O(\delta^3),
\end{equation}
where the $\Lambda_m \in \mathbb{C}$ satisfy
$\Lambda_m^*=\Lambda_{-m}$, with ${}^*$ denoting complex conjugate.
The multi-dimensional parameter $\Lambda$ is therefore 
infinite dimensional: $\Lambda=\{\Lambda_m \in\Comp : m \in\Zahl\}$.
Area preservation at $O(\delta^2)$ requires that $\Lambda_0 = 0$ and
the introduction of the $O(\delta^2)$, $\sig$-independent terms.

% ---------------------------------------------------------------------------

\subsection{Arbitrary axisymmetric flow}

The initial flow is taken to be axisymmetric, with vorticity
\begin{equation*}
   \omega_{\Lz}(r) = \frac{1}{r} \bigl(r \psi_\Lz'(r)\bigr)',
\end{equation*}
where the prime denotes differentiation with respect to $r$. 
For  this flow, \eqref{q:varphi1} reduces to
\begin{equation}\label{q:phi1e}
  \frac{\psi_\Lz'}{r} \Delta \partial_\sig \rho_1
	+ 2 \Bigl(\frac{\psi_\Lz'}{r}\Bigr)'
	  \Bigl(\partial^2_{r \sig} \rho_1
	  - \frac{1}{r} \partial_\sig \rho_1 \Bigr)
	  +  \frac{1}{r}\bigl(r \chi_1'\bigr)' = 0,
\end{equation}
with $\chi_1^{}$ a function of $r$ only.
The corresponding boundary condition is obtained from \eqref{q:bc1} in the form
\begin{equation}
   \partial_\sig \rho_1 = - \sum_m \Lambda _m \ex^{\i m \sig}
   \quad\textrm{at  } r=1.
\end{equation}
The solvability condition for \eqref{q:phi1e}, found by integration with
respect to $\sigma \in [0,2\pi]$, imposes that $(r \chi_1')'=0$;
boundedness of $\chi_1^{}$ then implies that $\chi_1^{}$ is a constant which we
can take equal to zero: $\chi_1^{}=0$.

The solution of \eqref{q:phi1e} for $\rho_1$ is then found as the
Fourier series
\begin{equation}\label{q:rho1fou}
   \rho_1(r,\sig) = \sum_{m} \Lambda_m {\rho}_{1,m}(r)\, \ex^{\i m\sig},
\end{equation}
with ${\rho}_{1,m}^*={\rho}_{1,-m}$.
Equation \eqref{q:phi1e} does not constrain the $m=0$ mode ${\rho}_{1,0}$;
this is the result of the gauge freedom for $\gL$.
A convenient choice is
\begin{equation} \label{q:mzero}
{\rho}_{1,0}=0.
\end{equation}
Introducing \eqref{q:rho1fou} into \eqref{q:phi1e} gives a second-order
equation for ${\rho}_{1,m}$, namely
\begin{equation}\label{q:phi1r}
  \psi_\Lz' \Bigl( {\rho}_{1,m}'' -\frac{1}{r} {\rho}_{1,m}'
	+ \frac{2-m^2}{r^2} {\rho}_{1,m}  \Bigr)
	+ 2 \psi_\Lz'' \Bigl( {\rho}_{1,m}'
	- \frac{1}{r} {\rho}_{1,m} \Bigr) = 0,
\end{equation}
with associated boundary condition
\begin{equation}\label{q:bc11}
  {\rho}_{1,m} = \frac{\i}{m}
  \quad\textrm{at  } r=1.
\end{equation}
There is a close connection between this equation and the Rayleigh equation
for the normal modes of axisymmetric flows \cite[e.g.,][]{drazin-reid:hs}:
\eqref{q:phi1r} can be recast as the Rayleigh equation for
zero-frequency modes, with $r^{-1} \psi_\Lz' {\rho}_{1,m}$ as the
unknown function. 
Of course the non-homogeneous boundary condition for ${\rho}_{1,m}$
differs from the homogeneous boundary condition usually considered for
the Rayleigh equation.
The connection is useful nevertheless: the absence of zero-frequency
normal modes that can be established from the Rayleigh equation when
$\psi_0' \not=0$ (as guaranteed by the hypothesis H2) implies the
existence of a unique solution to \eqref{q:phi1r}. 

We note that the solution for the $m=1$ mode, which describes a
rigid translation of the disc, is independent of $\psi_\Lz$ and given by
${\rho}_{1,\pm 1}=\i\Lambda_{\pm 1} r$.
Not surprisingly, this corresponds to a uniform
displacement field $\nabla^\bot \rho_1$. 

The vanishing of $\chi_1^{}$ indicates that the vorticity--streamfunction
relationship is unchanged at leading order in $\delta$.
This is a particularity of axisymmetric flows 
which makes it worthwhile to carry out the calculation to
$O(\delta^2)$ so as to demonstrate how a non-zero $\chi_2^{}$ is obtained;
this is described in Appendix~\ref{a:sec}. 

With $\rho_{1}$ determined by its Fourier series \eqref{q:rho1fou},
$\bar\phiform_1$ is given by
\begin{equation}\label{q:axiphiform1}
   \bar\phiform_{1}
   = \sum_{m} {\rho}_{1,m}(r) \ex^{\i m \sigma} \df \Lambda_{m}.
\end{equation}
Because ${\rho}_{1,0}=0$ and $1-\Ppsiz$ is simply the average along circles,
\eqref{q:noname1} indicates that $(1-\Ppsiz) \Psib1_1=0$.
Equations \eqref{q:phibs} then imply that
$\phinatb_{1}=\bar\phiform_{1}$.
In other words, our choice \eqref{q:mzero} provides the leading-order
connection with its natural choice of gauge which corresponds to
vanishing average along the circles $r=\textrm{const}$.
Expanding $\Psib1_1$ in Fourier series as
\begin{equation}
   \Psib1_{1}=\sum_{m} \Psib1_{1,m}(r)\, \ex^{\i m \sigma} \df \Lambda_{m},
\end{equation}
\eqref{q:noname1} is reduced to the set of ordinary differential equations
\begin{equation}\label{q:anyname}
  \psi_{\Lz}'\Bigl[\frac1r\bigl(r\Psib1_{1,m}\phantom{|}'\bigr)'
		- \frac{m^2}{r}\Psib1_{1,m}\Bigr]
	- \omega_{\Lz}' \Psib1_{1,m}
	=  - \omega_{\Lz}' {\rho}_{1,m}
\end{equation}
with $\Psib1_{1,0}=0$. 
The associated boundary conditions are found from \eqref{q:Psiobc} as
\begin{equation}
   \Psib1_{1,m} = \frac{\i}{m} \quad \textrm{at} \ \ r=1.
\end{equation}
Solving \eqref{q:anyname} gives the first-order correction $\Psio_{1}$
to the Eulerian flow to leading order in $\delta$.

As discussed at the end of \S5, the computation of the geometric angle to
leading order requires not only $\phinatb_{1}$ but also $\phinatb_{2}$ or, 
to minimise computations, $(1-\Ppsiz) \Psib1_{2}$.
This is deduced from \eqref{q:avpsiob}  which reduces to the ordinary
differential equation
\begin{equation}\label{q:phiform2*}
  \frac{1}{r} \frac{{\rm d}}{{\rm d}r} \Bigl\{ r  \frac{{\rm d}}{{\rm d}r}
	(1-\Ppsiz) \Psib1_{2} \Bigr\}
  = (1-\Ppsiz) \bigl\{ \Delta [\rho_1,\Psib1_{1}]
	- [\rho_1,\Delta \Psib1_{1}] \bigr\}.
\end{equation}
Solving this equation leads to an expression for $(1 - \Ppsiz) \Psib1_{2}$.
Taking the differential yields the first term of the curvature $\kappa$
in \eqref{q:kapphinat} as
\begin{equation} \label{q:dfhat}
   \langle \df \phinatb \rangle = (1 - \Ppsiz) \df \Psib1_{2} + O(\delta^{3}).
\end{equation}
Note that since $\rho_1$ is linear in $\Lambda_m$ and $\Psib1_1$ is
$\Lambda$-independent, $(1-\Ppsiz) \Psib1_{2}$ is linear in $\Lambda_m$;
furthermore, because the averaging $1-\Ppsiz$ along circles eliminates
all products in the right-hand side of  \eqref{q:phiform2*} except for
those of complex-conjugate Fourier modes, $(1-\Ppsiz) \Psib1_{2}$ is
a linear combination of terms of the type $\Lambda_m \df \Lambda_m^*$.
Therefore, $\langle \df \phinatb \rangle$ is given by a
$\Lambda$-independent linear combination of the two-forms
$\df\Lambda_m \wedge \df\Lambda^*_m$.

The second term in \eqref{q:kapphinat} is also $O(\delta^{2})$ and
is readily computed from \eqref{q:axiphiform1}.
Averaging along circles gives it the same form as that of
$\langle\df \phinatb \rangle$.
This leads to the geometric angle in the form
\begin{equation}\label{q:geoaxi}
   \geop = \delta^2 \sum_{m>0} f_m(r) \mathcal{A}_m + O(\delta^3),
\end{equation}
for some functions $f_m(r)$.
Here we have defined 
\[
   \mathcal{A}_m = -\frac{\i}{2} \int_{\mathcal{D}_\Lambda}
	\df\Lambda_m \wedge \df\Lambda^*_m
\]
which be recognised as (minus) the oriented area enclosed by the path
described by $\Lambda_m$ in the complex plane.
(A positive $\mathcal{A}_m$ is associated with a rotation of the fluid
domain in the positive sense.)
Unsurprisingly, at leading order, the geometric angle is the sum of
separate contributions of each Fourier mode of the boundary deformation.

% ---------------------------------------------------------------------------

\subsection{An example: flow with power-law radial dependence}

\begin{figure}
\begin{center}
\begin{tabular}{ccc}
$m=0$ & & $m=2$ \\
\epsfig{file=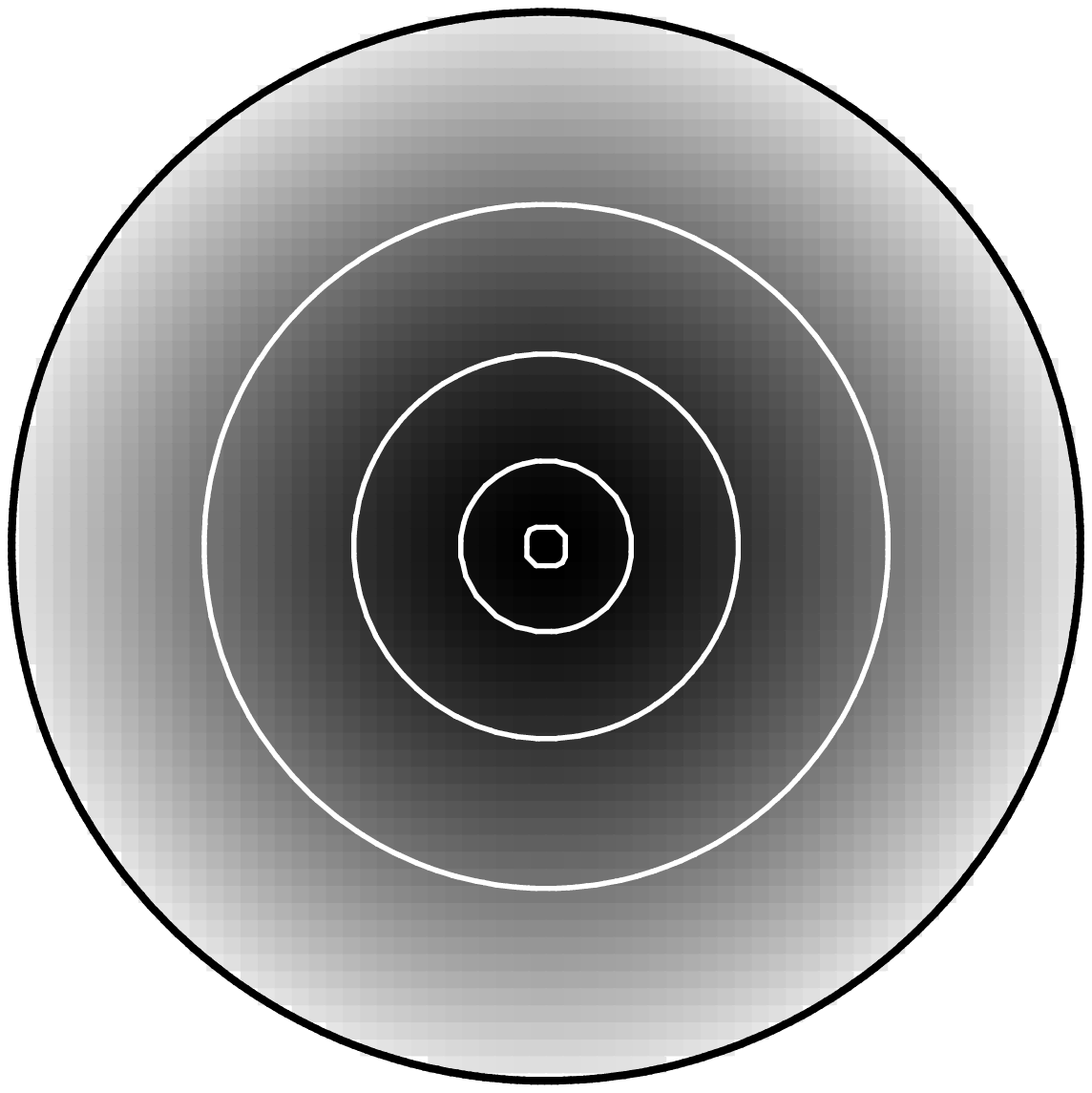,height=4.5cm} & \hspace{1cm} &
\epsfig{file=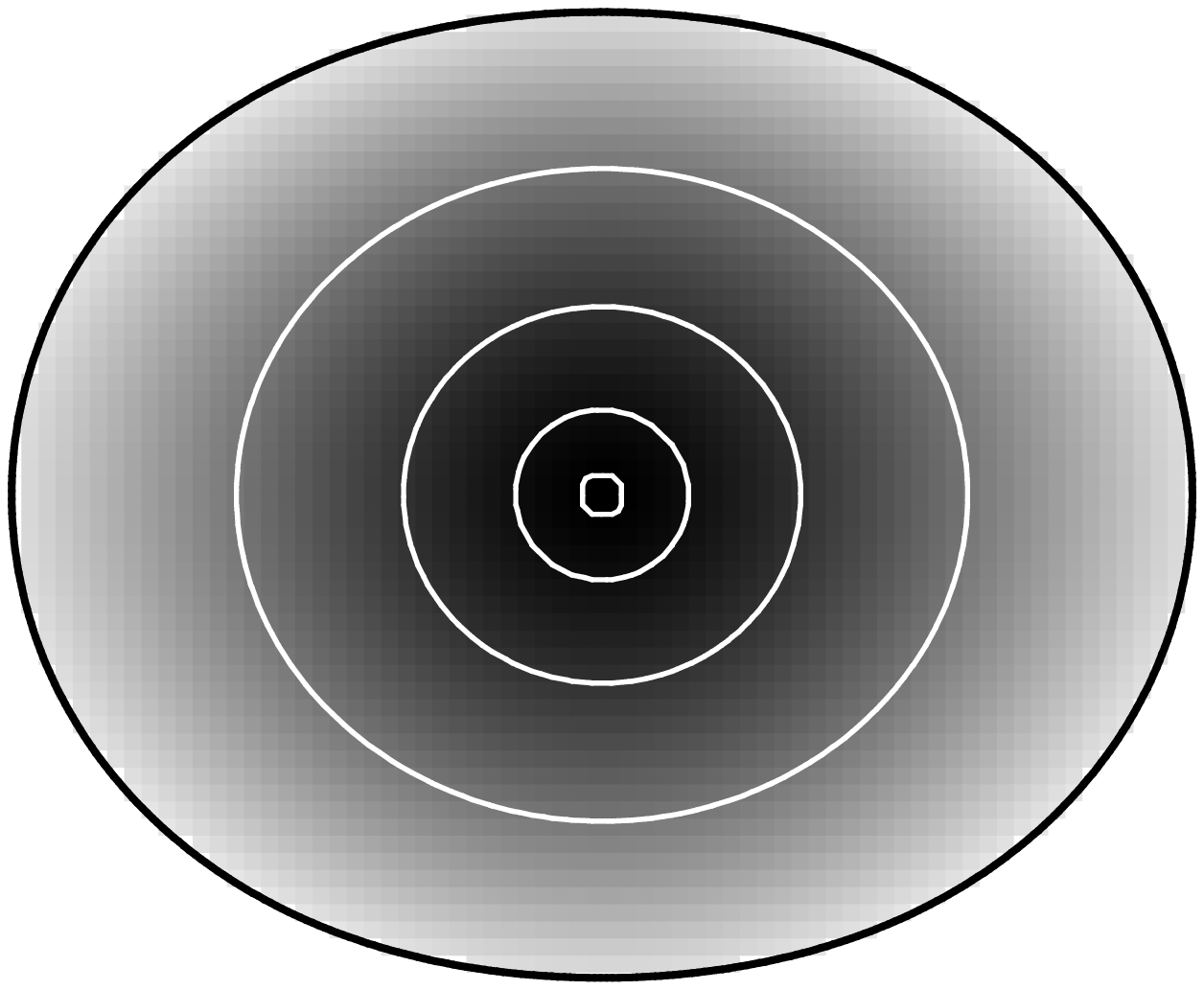,height=4.5cm} \\ 
$m=3$ & & $m=4$ \\
 \epsfig{file=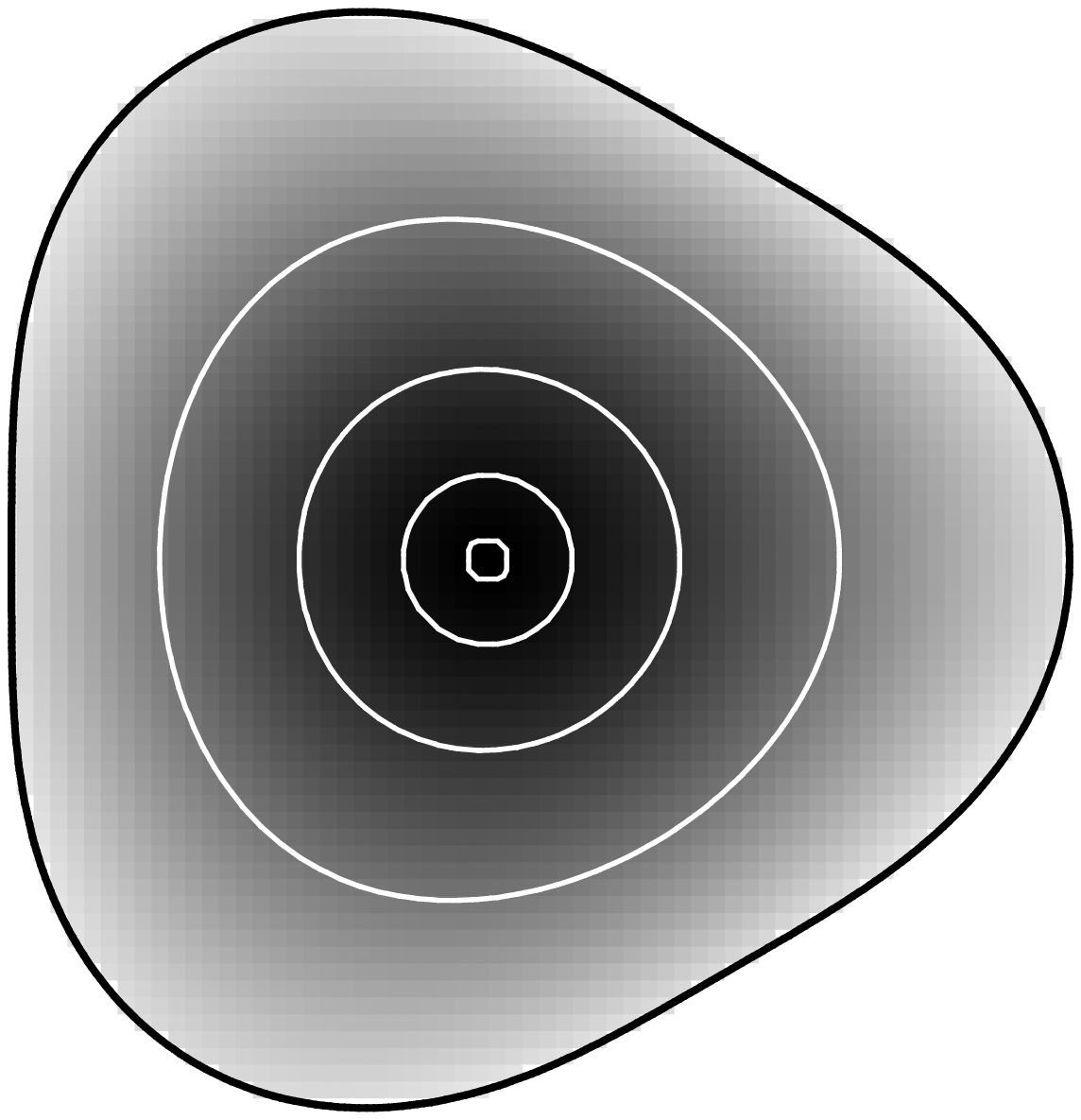,height=4.5cm} & &
 \epsfig{file=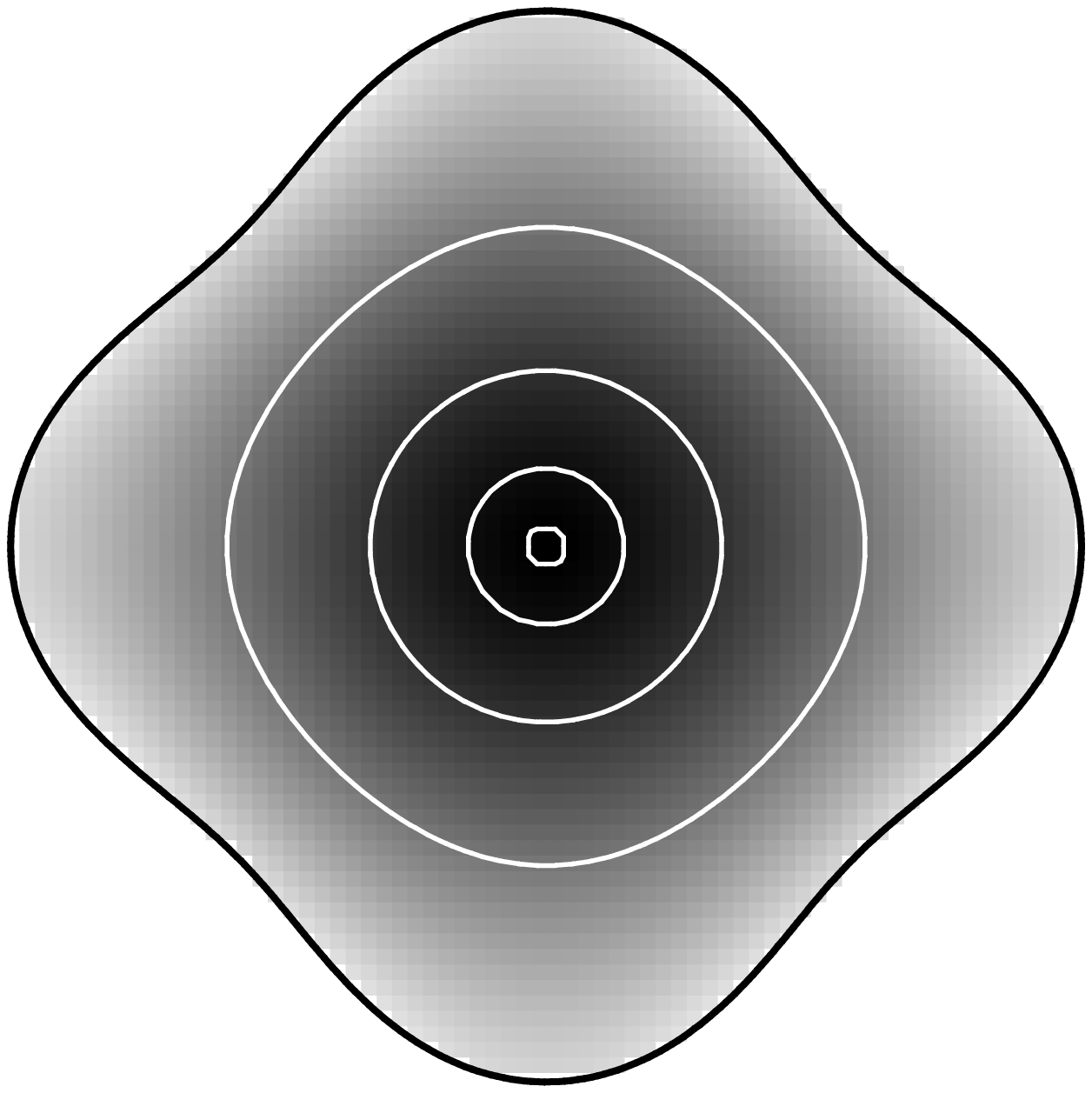,height=4.5cm} \\
$m=5$ & & $m=6$ \\
 \epsfig{file=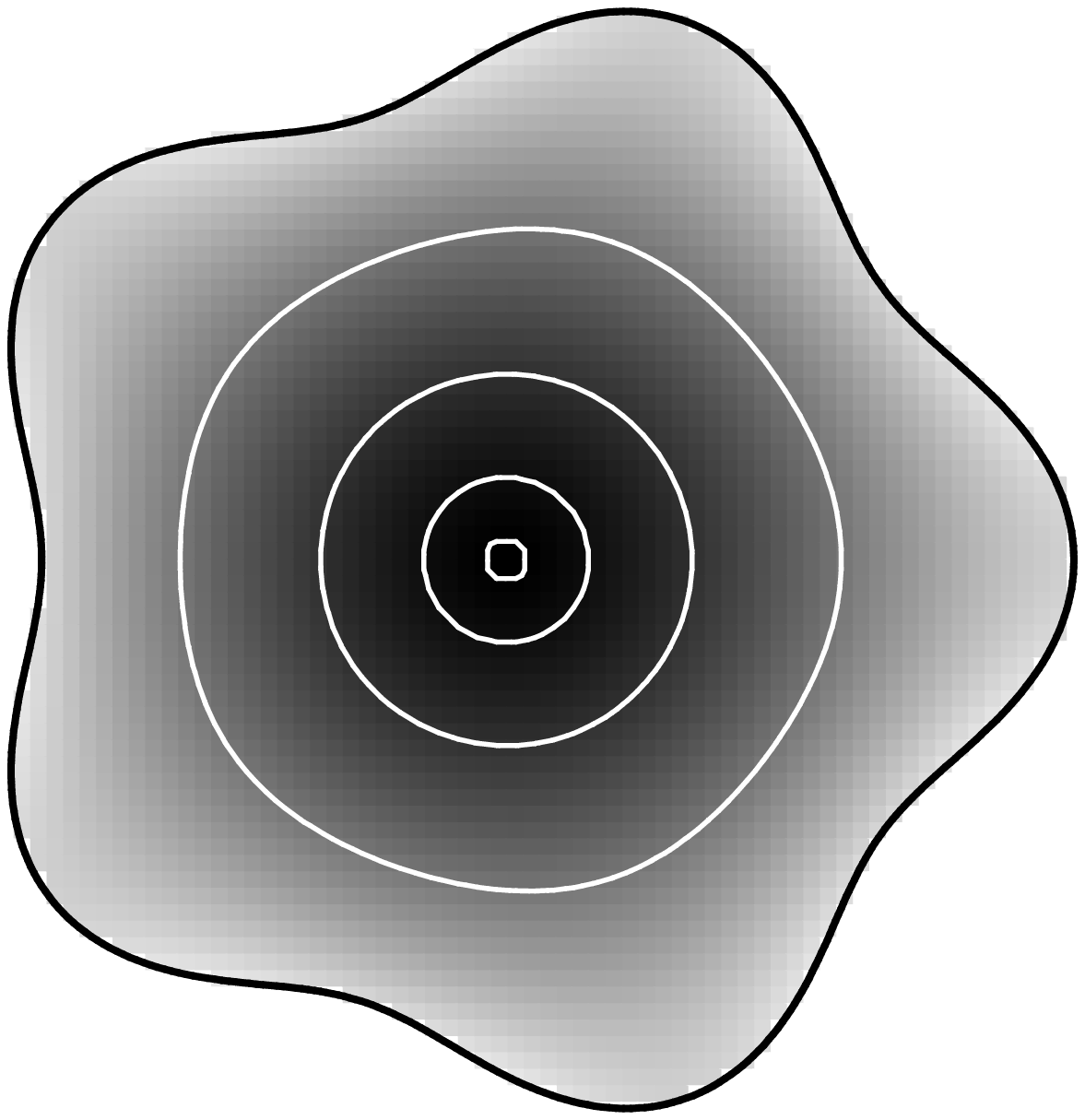,height=4.5cm} & &
 \epsfig{file=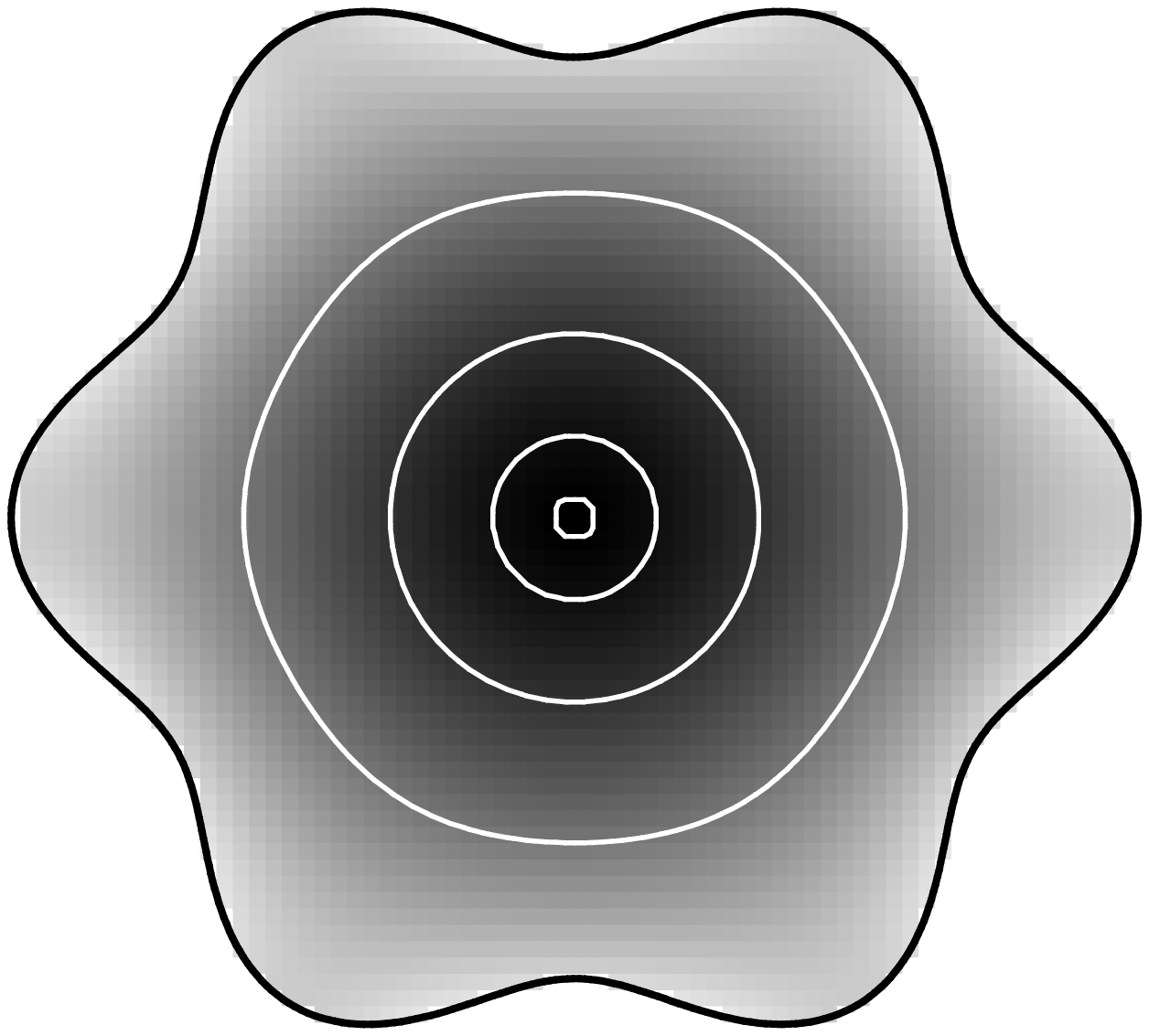,height=4.5cm} 
\end{tabular}
\end{center}
\caption{%
Vorticity $\wL$ (grey scale) and streamfunction $\psiL$ (white lines)
of the leading-order (steady) flow in a slightly deformed disc.
The top left panel $m=0$ shows the undeformed flow, with $\psi_\Lz=r^{1/2}$;
the other panels show the flows obtained when deforming the disc by
a single Fourier mode $m$ according to \eqref{q:disc} with
$\delta |\Lambda_m|=0.05$.}\label{f:disc}
\end{figure}
 
As a simple example of an axisymmetric flow, consider the streamfunction
\begin{equation}\label{q:axiexa}
   \psi_\Lz(r)=A r^\alpha
   \qquad\textrm{with} \ \  0<\alpha<2,
\end{equation}
for which \eqref{q:phi1r}--\eqref{q:bc11} can be solved explicitly, leading to
\begin{equation}\label{q:phi1ex}
   \bar\phiform_{1,m} = {\rho}_{1,m} = \frac{\i r^{\beta_m}}{m},
\end{equation}
where
\begin{equation*}
  \beta_m = \alpha_m-\alpha+2
  \quad\textrm{and}\quad
  \alpha_m = \sqrt{m^2+\alpha^2-2\alpha}.
\end{equation*}
The leading-order vorticity in the deformed domain is then found to be
\begin{equation}
   \wL(r,\sigma) = \omega_\Lz(r) - \delta\, \omega_\Lz'(r) \sum_{m} \Lambda_m
	r^{\beta_m-1} \ex^{\i m \sig} + O(\delta^2),
\end{equation}
with a similar expression for $\psiL$.
Figure \ref{f:disc} shows these approximations to $\wL$ and $\psiL$ in
domains deformed by a single Fourier mode $m$, with $m$ ranging from $2$
to $6$, in the case $\alpha=1/2$.
 
Equation \eqref{q:anyname} can also be solved explicitly, with the result
\begin{equation}
  \Psib1_{1,m} =  \frac{\i}{m}
	\bigl[ \gamma_m r^{\beta_m}+(1-\gamma_m) r^{\alpha_m}\bigr]
  \qquad\textrm{where }
  \gamma_m=\frac{\alpha}{\alpha_m+\beta_m}.
\end{equation}
Introducing into \eqref{q:phiform2*} gives
\begin{align}
  (1-\Ppsiz) \Psib1_{2} 
  &= -\i \sum_m \frac{1}{m} \Lambda_m \df {\Lambda}_m^*
	\bigl[ \gamma_m F(\beta_m,\beta_m) r^{2\beta_m-2}\notag\\
	  &\hbox to90pt{} {}+ (1-\gamma_m) F(\alpha_m,\beta_m)
		r^{\alpha_m+\beta_m-2}\bigr],
\end{align}
where we have defined
\begin{equation}\label{q:FG}
  F(\alpha_m,\beta_m)
  := \frac{E(\alpha_m,\beta_m)}{\alpha_m+\beta_m-2}
  := \frac{2\alpha_m\beta_m + \beta_m^2-2\alpha_m -2\beta_m+m^2}
	{\alpha_m+\beta_m-2}.
\end{equation}

From this and \eqref{q:dfhat} we deduce the first component of $\kappa$, namely
\begin{equation}\label{q:dphi2}\begin{aligned}
   \langle{\df\phinatb}\rangle &= - 2 \i \delta^2
	\sum_{m>0} \frac1m \bigl[\gamma_m F(\beta_m,\beta_m) r^{2\beta_m-2}
	\\
	&\hbox to2cm{}+ (1-\gamma_m) F(\alpha_m,\beta_m) r^{\alpha_m+\beta_m-2}\bigr] \;
		\df\Lambda_m \wedge\df \Lambda_m^* + O(\delta^3).
\end{aligned}\end{equation}
Using \eqref{q:phi1ex}, the second component of $\kappa$ is found directly
to be
\begin{equation}\label{q:phiwphi2}\begin{aligned}
   \langle{[\phinatb \wedge \phinatb ]}\rangle
	&= \frac{\delta^2}{2\pi} \int_0^{2\pi} 
		[\bar \phiform_1 \wedge \bar \phiform_1 ] \,
		\mathrm{d}\sigma + O(\delta^3) \\
	&= - 4  \i \delta^2 \sum_{m>0} \frac{\beta_m}{m}
		r^{2\beta_m-2} \, \df \Lambda_m \wedge \df\Lambda^*_m
		+ O(\delta^3).
\end{aligned}\end{equation}

\begin{figure}
\begin{center}
\epsfig{file=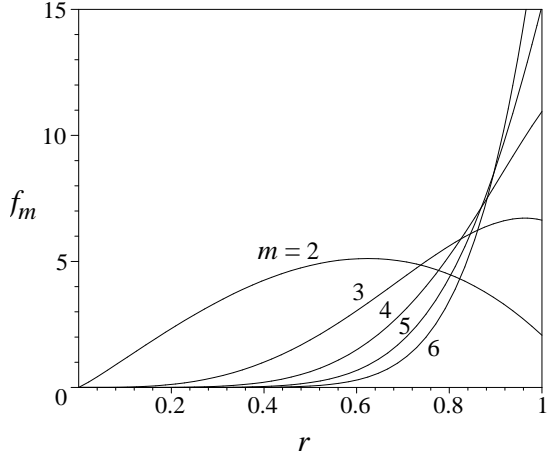,height=6cm}
\caption{%
Functions $f_m(r)$ giving the geometric angle in \eqref{q:geoaxi}
for the axisymmetric flow with streamfunction $\psi_0=r^{1/2}$ in a disc
deformed by a Fourier mode $m$.}\label{angle}
\end{center}
\end{figure}

Combining these results with \eqref{q:geop3}, and noting that the
action--angle variables in the undeformed domain are simply
$(I,\theta)=(r^2/2,\sigma)$, leads to the geometric angle in
the form \eqref{q:geoaxi} with
\begin{equation} \label{q:fm}
   f_m(r) = \frac{4}{m} \left( p_m r^{2\beta_m-4}
     + q_m r^{\alpha_m+\beta_m-4}\right),
\end{equation}
where
\begin{equation}\label{q:pmqm}
    p_m = \gamma_m E(\beta_m,\beta_m)-2\beta_m(\beta_m-1) 
    \quad\textrm{and}\quad
    q_m = (1-\gamma_m)E(\alpha_m,\beta_m).
\end{equation}
Figure~\ref{angle} shows the functions $f_m(r)$ for $m=2,3,\cdots,6$ in the
case $\alpha=1/2$. 

A spot check for our results is provided by the limit $\alpha\to2$,
corresponding to a flow with a uniform vorticity $4A$.
Assuming that $\Lambda_m=0$ if $m \not= \pm 2$ and  that $\Lambda_{\pm2}$ trace
the unit circle in the complex plane, the domain deformation is simply
the rotation of a small-eccentricity ellipse with semi-axes $1+2\delta$
and $1-2\delta$.
There is an exact analytic solution for such a uniform-vorticity flow
in a rotating ellipse; both the direct use of this solution and
\eqref{q:geoaxi} yield the same $r$-independent value for the
leading-order geometric angle that appears for one full rotation of
the ellipse, namely $16 \pi \delta^2$.
See Appendix~\ref{a:rot} for details.

% ==========================================================================

\section{Discussion}

In this paper, we have used differential geometry mostly as a notational
tool, its main application being to make explicit the linear dependence of
various quantities on $\Lambdadot$.
It is nonetheless clear that the objects we are dealing with can be given
an interpretation in a more abstract geometric setting and that the
problem may be placed in the framework of geometric mechanics.
This would be a valuable undertaking, but one which is probably quite
involved due to the infinite-dimensional nature of the problem,
and which is certainly beyond the scope of the present paper.
However, we now discuss briefly and informally the geometric context of
our results in order to elucidate the meanings of the connection
one-form $\phinat$, the form of its curvature $\kappa$, etc.

Central to our development is the group $G=\mathrm{SDiff}(\RReal^2)$
of area-preserving diffeomorphisms of the plane.
Taking the initial domain $D_\Lambdaz$ as a reference domain,
the subgroup $H\subset G$ which maps $D_0$ to itself, viz.,
\begin{equation}
   H = \{ g \in G : g D_\Lambdaz = D_\Lambdaz\},
\end{equation} 
is of particular importance.
In terms of $G$ and $H$, our parameter space $\Lspace$, the space of
all possible shapes $\Lambda$ of the domain $D_\Lambda$, can be realised
as the quotient $G/H$.
Indeed, each right coset of $H$ in $G$ contains all diffeomorphisms
mapping $D_0$ to $D_\Lambda$ for a particular $\Lambda$, since any
two such diffeomorphisms $g$ and $g'$ are related by $g'=g h$
for some $h \in H$.
We can therefore identify $\Lspace$ with $G/H$.
Another important subgroup is 
\begin{equation}
   H_0 = \{ h \in H : \omega_\Lambdaz \circ h = \omega_\Lambdaz \}  
\end{equation}
containing all area-preserving diffeomorphisms in $D_\Lambdaz$ which
leave the initial vorticity distribution invariant.
A rearrangement $\omega_\Lambdaz \circ g $ of the initial vorticity
$\omega_\Lambdaz$ can be identified with an element of $G/H_0$.

Armed with this setup, we can interpret our results geometrically.
The key point is to regard $G$ and $G/H_0$ as principal bundles, both
with $G/H\simeq \Lspace$ as base manifold.
Finding the (leading-order) Eulerian flow for each domain shape $\Lambda$
then amounts to finding a lift from $G/H$ to $G/H_0$;
finding the Lagrangian particle position amounts to finding a lift
from $G/H$ to $G$.\footnote{%
In our geometric description the interior and exterior of $D_\Lambda$
are treated on the same footing; this can be done because our formulas,
with suitable boundary conditions as $|\xb| \to \infty$, would also
apply to a fluid flow outside $D_\Lambda$.}
Proposition P1, stating that the leading-order Eulerian flow depends
only on the domain shape, says that the lift from $G/H$ to $G/H_0$
is path-independent;
in other words, it defines a section of $G/H_0$.\footnote{%
We stress again the local nature of P1: globally, there are many possible
steady flows for a given vorticity distribution and a given domain $D_\Lambda$.
Geometrically, this implies that the section of $G/H_0$ is multivalued.}
In contrast, the lift from $G/H$ to $G$, which gives an approximation to
the particle position, depends on the path in $G/H$ and in fact on the
speed with which the path is traced.
There is, however, a contribution that is independent of speed;
for cyclic domain deformations, it is quantified by the geometric
angle given in P2.

It is worth commenting on the meaning of the one-form $\phinat$ that
appears in the geometric angle.
One way of defining a lift in a principal bundle is by means of a
vector-valued one-form, i.e.\ a linear map from $T(G/H)$ to $TG$,
describing the vertical (along-fibre) displacement associated with
any given displacement on the base manifold $G/H$.
Such a form can be recognised as a connection form.
In our context, $TG$ is the space of divergence-free vector fields over
$\RReal^2$, which can be identified through the use of a streamfunction
with the space of real-valued functions $\CFun(\RReal^2)$.
Thus, a lift can be defined by a connection one-form over $G/H\simeq \Lspace$
with values in $\CFun(\RReal^2)$.
This is precisely the interpretation we give to $\phinat$.
With the geometric interpretation of $\phinat$, the subsequent results are
clear: \eqref{q:kappa} is the standard expression for the curvature of
$\phinat$, the geometric angle \eqref{q:geop4} is given by the holonomy of
$\phinat$, and the standard conclusion about the geometrical angle in
finite-dimensional systems is recovered.

Throughout the paper, we emphasize that our results are local in nature:
the prediction of an instantaneously steady Eulerian flow, for
instance, holds only if the domain deformation is such that H1 and H2 are
always satisfied.
This is necessarily the case if they are satisfied initially and the domain
deformation is sufficiently small, but it may well continue to hold for
larger deformations.
It is nonetheless interesting to speculate about the dynamics when either
H1 or H2 fails in the course of the evolution.
If H1 fails, the flow ceases to be Arnold stable and likely becomes
spectrally unstable.
We can then expect the flow to become highly unsteady and, in the absence
of dissipative mechanisms, remain so regardless of subsequent deformations
of the domain.

The failure of H2, on the other hand, corresponds to the appearance of
streamlines for which the orbiting period of particles becomes large.
When the period becomes comparable with the time scale of the domain
deformation, our asymptotic approach clearly breaks down.
This can happen when the flow is driven by the domain deformation towards a
change in topology, with the creation of hyperbolic stagnation points and
separatrices.
How the flow evolves in this situation is unclear, but some understanding
could be gained by investigating the problem where the initial steady
flow $\psi_\Lambdaz$ has a hyperbolic stagnation point (and satisfies
H1---the Kelvin--Stuart vortex in \cite{holm-al:stab:85} is just one example).
This problem can be viewed as a generalisation of the classical
critical-layer problem for parallel shear flows (e.g.,
\cite{stewartson:78,maslowe:86} and references therein).
In this generalisation, the separatrix plays the role of the zero-velocity
critical line (along which H2 is obviously violated); by analogy, it can be
expected to be also surrounded by a narrow critical layer where complicated
nonlinear dynamics occurs.
We plan to investigate this problem in future work.

We conclude this paper by remarking on the possible extension of our
results to flows in three dimensions.
In three dimensions, the dynamics of an inviscid and incompressible fluid
is determined, as in two dimensions, by a form of conservation of
vorticity, although in this case it is as a vector that the vorticity
is transported \cite[e.g.,][]{arnold-khesin:tmh}.
This suggests that our approach for the determination of the leading-order
Eulerian flow in deforming domains can be adapted to the three-dimensional
setting.
The technical conditions for the well-posedness of the
equations for $\gL$ are however likely to be significantly more
complicated than in two dimensions.

The evolution of the fluid-particle position seems, at first sight, to pose
a very different problem in three than in two dimensions, since the
velocity field is divergence-free and not Hamiltonian.
However, particle trajectories for (non-Beltrami) steady solutions of the
Euler equations are known to be integrable \cite[]{arnold:65,arnold:66}
because they are confined to surfaces of constant Bernoulli function (Lamb
surfaces).
There is, therefore, a simple characterisation of the fluid-particle
positions in steady flows, analogous to the action--angle characterisation
in two dimensions.
This could be used for slowly time-dependent flows to quantify the effects
of a cyclic boundary deformation as was done to obtain the geometric angle
in this paper.

A general difficulty with three-dimensional flows, however, is the absence
of general stability results similar to those obtained by the
energy--Casimir method \cite[]{holm-al:stab:85}.
Instabilities cannot therefore be excluded (on the contrary, they are the
rule rather than the exception), and the effect of their competition with
the slow evolution of the leading-order flow would need to be assessed
carefully.

% ===========================================================================

% \begin{acknowledgements}
\medskip\noindent{\bf Acknowledgements.}\ \
This research was supported by an EPSRC research grant, a William Gordon
Seggie Brown Fellowship (DW) and a NERC Advanced Fellowship (JV).
Additional support was provided by the EPSRC network ``Geometrical method
in geophysical fluid dynamics''.
We thank T.~N.~Bailey, Y.~Brenier and L.~Butler for helpful discussions.
% \end{acknowledgements}

% ===========================================================================

\appendix

\section{Properties of $\phiform$}\label{a:details}

In this Appendix, we give details of the derivations
of three useful expressions. 

\medskip\noindent{\bf 1.}
We first derive \eqref{q:integr}.
Applying $\df$ to \eqref{q:dfwL} gives
\begin{equation}\label{q:precurv1}\begin{aligned}
   0 = \df^2\wL &= -\df[\phiform,\wL]\\
	&= -[\df\phiform,\wL] + [\phiform\wedge\df\wL]\\
	&= -[\df\phiform,\wL] - [\phiform\wedge[\phiform,\wL]],
\end{aligned}\end{equation}
where the bracket $[\cdot\wedge\cdot]$ is defined in \eqref{q:wbra}.
Now
\begin{equation}\begin{aligned}
   -[\phiform\wedge[\phiform,\wL]]
	&= -[\phiform_m,[\phiform_n,\wL]]\; \dLwdL{m}{n}\\
	&= \bigl\{ [\phiform_n,[\wL,\phiform_m]]
	   + [\wL,[\phiform_m,\phiform_n]] \bigr\}\; \dLwdL{m}{n}\\
	&= [\phiform\wedge[\phiform,\wL]] + [\wL,[\phiform\wedge\phiform]]
\end{aligned}\end{equation}
(we have used Jacobi's identity, $[f,[g,h]]+[g,[h,f]]+[h,[f,g]]=0$ for any
three functions $f$, $g$ and $h$, to arrive at the second equality);
therefore
\begin{equation}\label{q:Phi-id0}
   -[\phiform\wedge[\phiform,\wL]]=\sfrac12[\wL,[\phiform\wedge\phiform]]
\end{equation}
and from \eqref{q:precurv1} we find
\begin{align}
   &[\df\phiform + \sfrac12 [\phiform\wedge\phiform], \wL] = 0  \notag\\
   \Rightarrow\qquad
   &\df\phiform + \sfrac12 [\phiform\wedge\phiform]
	= w\circ\wL \label{q:curv}
\end{align}
for an arbitrary function-valued two-form $w \circ \wL$.

% ---------------------------------------------------------------------------

\medskip\noindent {\bf 2.}
Next we show that \eqref{q:curv} can be derived directly
from \eqref{q:pde0}.
This provides a consistency check for our developments.
Let us first compute
\begin{equation}\begin{aligned}
   \df\,[\phiform,\psiL] &= \df\bigl\{ [\phiform_m,\psiL]\,\df\Lambda_m \bigr\}\\
	&= [\df\phiform_m,\psiL]\wedge\df\Lambda_m
	  + [\phiform_m,\df\psiL]\wedge\df\Lambda_m\\
	&= [\df\phiform,\psiL] - [\phiform\wedge\df\psiL]
\end{aligned}\end{equation}
and
\begin{equation}\begin{aligned}
   \df(\df\GL\circ\wL) &= -(\df\GL'\circ\wL)\wedge\df\wL
	= (\df\GL'\circ\wL)\wedge[\Phi,\wL]\\
	&= [\phiform\wedge\df\GL],
\end{aligned}\end{equation}
where the last equality can be verified by computation in coordinates.

Writing \eqref{q:pde0} as
\begin{equation}
   \lapl\,[\phiform,\psiL] - [\phiform,\wL] - \lapl(\df\GL\circ\wL) = 0
\end{equation}
and taking $\df$, we find
\begin{align}
  0	&=
   \lapl\,\df[\phiform,\psiL] - \df[\phiform,\wL] - \lapl\df(\df\GL\circ\wL)
		\notag\\
	&= \lapl[\df\phiform,\psiL] - \lapl[\phiform\wedge\df\psiL]
	   - [\df\phiform,\wL] + [\phiform\wedge\df\wL]
	   + \lapl[\phiform\wedge\df\GL] \notag\\
	&= (\lapl-\FL')\,[\df\phiform,\psiL]
	  + \lapl\,[\phiform\wedge(\df\GL-\df\psiL)] + [\phiform\wedge\df\wL]
		\notag \\
	&= (\lapl-\FL')\,[\df\phiform +\sfrac12[\phiform\wedge\phiform],\psiL].
\end{align}
A couple of identities have been used to arrive at the last equation.
The first one is
\begin{equation}\label{q:Phi-id2}
   [\phiform\wedge[\phiform,\psiL]]
	= \sfrac12\,[[\phiform\wedge\phiform],\psiL],
\end{equation}
which is proved in the same way as \eqref{q:Phi-id0}.
The second identity is
\begin{equation}\label{q:Phi-id1}
\begin{aligned}
\left[\phiform \wedge \df \wL\right] &= [\phiform\wedge[\phiform,\FL\circ\psiL]]
	= [\phiform\wedge\FL'[\phiform,\psiL]]\\
	&= \FL'[\phiform\wedge[\phiform,\psiL]]
	   - [\phiform,\FL']\wedge[\phiform,\psiL]
	= \FL'[\phiform\wedge[\phiform,\psiL]],
\end{aligned}\end{equation}
where we have used $[\phiform,\FL'\circ\psiL]\wedge[\phiform,\psiL]
= \FL''\,[\phiform,\psiL]\wedge[\phiform,\psiL] = 0$ for the last equality.

The desired result \eqref{q:curv} is recovered by noting that
the operator $(\lapl-\FL')$ is invertible by hypothesis and that
$[\df\phiform +\sfrac12[\phiform\wedge\phiform],\psiL]=0$ on $\dy D_\Lambda$.
The latter can be verified by differentiating \eqref{q:dfwL} and evaluating it
on $\dy D_\Lambda$.

% ---------------------------------------------------------------------------

\medskip\noindent
{\bf 3.} Finally, we establish the formula
\begin{equation}\label{q:prop3}
   \df\phiformb - \sfrac12[\phiformb \wedge\phiformb]
	= \df\phiform + \sfrac12[\phiform \wedge\phiform].
\end{equation}
Its application to the natural connection $\phinat$ shows that $\kappa$
is independent of $\theta$.
Our proof starts by noticing that
$[\phiformb\wedge\phiformb]=[\phiform\wedge\phiform]$
because the transformation to action--angle variables is canonical.
Thus \eqref{q:prop3} is equivalent to
\begin{equation}
   \df\phiformb = \df\phiform + [\phiform \wedge\phiform].
\end{equation}
This is established by direct computation as follows
\begin{align}
   \df \phiformb(I,\theta;\Lambda)
	&= \df \phiform(\Xb(I,\theta;\Lambda);\Lambda) \notag\\
	&= \df \phiform(\xb;\Lambda)
	+ \sum_{n}\, \bigl(\df \Xb \cdot \nabla \phiform_n\bigr)\big|_{(\xb,\Lambda)} \,\df\Lambda_n \notag \\
	&= \df \phiform(\xb;\Lambda)
	+ \sum_{m,n}\, [\phiform_m,\phiform_n]\big|_{(\xb,\Lambda)} \;\dLwdL{m}{n} \notag\\
	&= \df \phiform(\xb;\Lambda) + [\phiform \wedge\phiform]\big|_{(\xb,\Lambda)}.
\end{align}

% ===========================================================================

\section{Solution of $(\lapl-\FL'\PpsiL)\,u=f$}\label{a:lin}

Here we show that the problem
\begin{equation}\label{q:DP}\begin{aligned}
   &(\Delta - \FL'\,\PpsiL)\,\eta = f\\
   &\eta = g \qquad\textrm{on }\dy D_\Lambda
\end{aligned}\end{equation}
has a unique solution $\eta$ when $\FL'  > - \cpoi$ everywhere in $D_\Lambda$ as follows from the hypothesis H1.

\medskip
We start with an identity.
Let $u$ and $v$ be such that $\PpsiL u=0$ and $\PpsiL v=v$.
We have
\begin{equation}\label{q:Porth}
   \int_D u\,v \dxdy
	= \int \Bigl\{ \oint u\,v \ds \Bigr\}\; {\rm d}\psiL
	= \int u\, \Bigl\{ \oint v \ds \Bigr\}\; {\rm d}\psiL
	= 0.
\end{equation}
From this it follows that the projection $\PpsiL$ is orthogonal
in $L^2(D_\Lambda)$,
in the sense that for any (sufficiently smooth) function $w$
\begin{equation}\label{q:Porth2}\begin{aligned}
   \int_D |w|^2 \dxdy
	&= \int_D \bigl\{ |\PpsiL w|^2 + 2 (\PpsiL w)[(1-\PpsiL) w]
		+ |(1-\PpsiL)w|^2 \bigr\} \dxdy\\
	&= \int_D \bigl\{ |\PpsiL w|^2 + |(1-\PpsiL)w|^2 \bigr\} \dxdy.
\end{aligned}\end{equation}

Using \eqref{q:Porth}, we find that the operator $(\Delta-\FL'\,\PpsiL)$
is self-adjoint for any functions $u$ and $v$ which vanish on $\dy D$,
\begin{equation}
   \int_D v\, (\Delta-\FL'\,\PpsiL)\, u \dxdy
	= \int_D u\, (\Delta-\FL'\,\PpsiL)\, v \dxdy.
\end{equation}
Moreover, $(\Delta-\FL'\,\PpsiL)$ is coercive under the hypothesis
$F'> - \cpoi$.
To show this, we first combine \eqref{q:Porth2} and Poincar{\'e} inequality
to obtain
\begin{equation}
   \int_D |\PpsiL u|^2 \dxdy
	\le \int_D |u|^2 \dxdy
	\le \frac{1}{\cpoi} \int_D |\gb u|^2 \dxdy
\end{equation}
for any function $u$ vanishing on $\dy D$.
It is then clear that
\begin{equation}
   \int_D u\,(\Delta-\FL'\,\PpsiL)\, u \dxdy
	= -\int_D \bigl\{ |\gb u|^2 + \FL'(\PpsiL u)^2 \bigr\} \dxdy
	\le 0
\end{equation}
when $\FL' > \cpoi$ everywhere in $D$, with equality obtaining
only when $u=0$.

\medskip
Returning to the problem \eqref{q:DP}, we extend $g$ to
$\textrm{cl}\, D_\Lambda$ and let $\tilde\eta=\eta-g$.
The problem thus becomes
\begin{equation}\label{q:hatvfi}
   (\Delta - \FL'\PpsiL)\,\tilde\eta = -\Delta g + \FL'\PpsiL g
\end{equation}
with boundary conditions $\tilde\eta=0$ on $\dy D_\Lambda$.
We have shown that the operator $(\Delta-\FL'\PpsiL)$ on the left-hand
side is self-adjoint and its associated bilinear form is coercive,
so by the Lax--Milgram lemma (assuming compactness, etc.,
cf.~\cite{gilbarg-trudinger:epde}) a unique solution $\eta$ exists
for \eqref{q:DP}.

% ===========================================================================

\section{Second-order terms in nearly axisymmetric flows}\label{a:sec}

At order $O(\delta^2)$, $\rho_2$ is found from \eqref{q:varphi2} to satisfy
\begin{equation}\label{q:phi22}\begin{aligned}
    \frac{\psi_\Lz'}{r} \Delta \partial_\sig \rho_2
	&+ 2 \Bigl(\frac{\psi_\Lz'}{r}\Bigr)' \Bigl(\partial^2_{r \sig} \rho_2
	- \frac{1}{r} \partial_\sig \rho_2 \Bigr)
	+ \frac{2}{r}\left(r \chi_2'\right)'\\
	&= 2 [ \rho_1,\Delta[ \rho_1,\psi_\Lz]]
	- [\rho_1,[\rho_1,\Delta \psi_\Lz]]
	- \Delta [\rho_1,[\rho_1,\psi_\Lz]].
\end{aligned}\end{equation}
The boundary condition \eqref{q:bc2} can be written as
\begin{equation}\label{q:phi21}
  \partial_\sig \rho_2
  = \partial_\sig \left( \partial_r \rho_1 \partial_\sig \rho_1 \right)
  -  \left( \partial_\sig \rho_1 \right)^2 +  \sum_m |\Lambda_m|^2
  \quad\textrm{at } r=1,
\end{equation}
after some manipulations.
The interest of this form is that, when \eqref{q:bc11} is taken into
account, it is clearly consistent, with both sides having a vanishing
$\sig$-average.
A solvability condition for \eqref{q:phi22} is obtained by averaging
over $\sig$, leading to
\begin{equation*}
 (r \chi_2')' = -  \Bigl(
  \frac{\psi_\Lz'}{2\pi r} \int_0^{2\pi} \Bigl[ (\partial_{r\sig}^2\rho_1)^2
    + \frac{1}{r^2} (\partial_{\sig\sig}^2 \rho_1)^2
    - \frac{2}{r} \partial_\sig \rho_1 \partial^2_{r\sig} \rho_1 \Bigr]
  \, \mathrm{d} \sig \Bigr)'.
\end{equation*}
This equation determines $\chi_2^{}$ uniquely up to an irrelevant
arbitrary constant.
When it is satisfied, \eqref{q:phi22} can be solved for $\rho_2$,
yielding a solution in the form of a Fourier series
\begin{equation*}
  \rho_2(r,\sig)=\sum_{m} \hat{\rho}_{2,m}(r)\, \ex^{\i m \sig},
\end{equation*}
with $\hat{\rho}_{2,m}^*=\hat{\rho}_{2,-m}$ and $\hat{\rho}_{2,0}=0$.
The functions $\hat{\rho}_{2,m}$ satisfy an inhomogeneous version of
\eqref{q:phi1r} obtained from \eqref{q:phi21};
clearly, they are quadratic in the $\Lambda_m$.

% ===========================================================================

\section{Rotating ellipse}\label{a:rot}

Consider a fluid inside an ellipse with semi-axes $a$ and $b$ that is
rotating around its centre with a (possibly time-dependent) angular
velocity $\eps \dot{\lambda}(\tau)$. The equation of the ellipse is given
by
\[
   B(x,\tau) = \frac{{\hat{x}_1}^2}{a^2}+\frac{{\hat{x}_2}^2}{b^2} - 1 = 0,
\]
where
\begin{align*}
  \hat{x}_1 &= x \cos \lambda  + y \sin \lambda, \\
  \hat{x}_2 &= -x \sin \lambda + y \cos \lambda
\end{align*}
are Cartesian coordinates in a frame rotating with angular velocity
$\eps\dot{\lambda}$.
An exact solution for the fluid motion in such a rotating ellipse is
provided by the uniform-vorticity flow with streamfunction
\[
   \psi(\xb,t) =
	K \Bigl(\frac{{\hat{x}_1}^2}{a^2}+\frac{{\hat{x}_2}^2}{b^2}\Bigr)
	+ \frac{\eps \dot{\lambda}(a^2-b^2)}{2(a^2+b^2)}
		\left({\hat{x}_1}^2-{\hat{x}_2}^2\right),
\]
where $K$ is a constant.
In this streamfunction, which can be verified directly to satisfy
the boundary condition \eqref{q:bc0}, the first term can be identified
with $\psiz$, the second with $\psio$, and there are no higher-order terms
in $\eps$ \cite[see, e.g.,][p.\ 421, for a derivation]{jeff-jeff}.
Action--angle coordinates $(I,\theta)$ for this flow satisfy  
\begin{align*}
  \hat{x}_1 = \sqrt{2Ia/b} \,\cos\theta,\\
  \hat{x}_2 = \sqrt{2Ib/a} \,\sin\theta.
\end{align*}
In terms of these variables, the streamfunction, or Hamiltonian, becomes
\begin{equation} \label{q:Hell}
  \hat{H}(I,\theta) = \frac{2 K I}{ab} + \frac{\eps\dot{\lambda}I}{ab}
	\Bigl[\frac{a^2-b^2}{a^2+b^2} (a^2 \cos^2 \theta - b^2 \sin^2 \theta)
	- (a^2 \cos^2 \theta + b^2 \sin^2 \theta) \Bigr].
\end{equation}
In this expression, the last term between round brackets comes from
the time-dependence in the (canonical) transformation from $(x,y)$
to $(I,\theta)$; it simply corresponds to a rigid-body rotation with 
angular velocity $-\eps \dot{\lambda}$.
The geometric angle is derived by writing the evolution equation for
$\theta$, averaging, then integrating in time.
This gives
\begin{equation} \label{q:ellgeo}
   \geop = \frac{\Delta \lambda}{2 ab}
	\Bigl[ \frac{(a^2-b^2)^2}{a^2+b^2} - (a^2+b^2) \Bigr],
\end{equation}
where $\Delta\lambda$ is the total angle rotated by the ellipse. 

Note that the averaging is in fact unnecessary, since the Hamiltonian
$\hat{H}$ does not depend on $\theta$ as a simplification of
\eqref{q:Hell} indicates. We do not perform this simplification here,
however, in order to retain the two terms in \eqref{q:ellgeo}
separately. This facilitates the comparison with the general formalism
of \S4. The first term in \eqref{q:ellgeo} stems from the correction in
the streamfunction 
$\psio$ (which here corresponds to a potential flow), while the
second stems from the slow time-dependence of the leading-order flow;
thus these two terms can be identified with the two contributions
$\langle\df\phinath\rangle$ and 
$-\frac{1}{2} \langle{[\phinath\wedge\phinath]}\rangle$ in \eqref{q:geop3}.

We can use the exact formula \eqref{q:ellgeo} to verify the approximate
results for slightly deformed axisymmetric flows obtained in \S6.
A rotating ellipse of small eccentricity is represented by the deformed
disc \eqref{q:disc} with $\Lambda_{\pm2}(\tau)$ describing a unit circle
in the complex plane and all the other $\Lambda_m$ equal to zero. 
The corresponding semi-axes are then
\[
   a = 1 + 2 \delta + O(\delta^2)
   \quad\textrm{and}\quad
   b = 1 - 2\delta + O(\delta^2),
\]
with the $O(\delta^2)$ corrections ensuring that $ab=1$. 
Introducing this into \eqref{q:ellgeo} and considering a full rotation
$\Delta\lambda = 2\pi$ gives the geometric angle
\begin{equation}\label{q:ellgeo2}
   \geop = \pi \left[ 32 \delta^2 - (2 + 16 \delta^2) \right] + O(\delta^3)
	 = - 2 \pi + 16 \pi \delta^2 + O(\delta^3). 
\end{equation}
An equivalent result is obtained from the developments in~\S6.
Since the only non-zero parameter $\Lambda_m$ is $\Lambda_{\pm2}$,
\eqref{q:geoaxi} reduces to
\[
   \geop = \delta^2 f_2(r) \mathcal{A}_2 + O(\delta^3).
\]
The uniform-vorticity corresponds to the limit $\alpha \rightarrow 2$ in  
\eqref{q:axiexa}, so that $\alpha_2=\beta_2=2$, $\gamma_2=1/2$, and
$f_2(r)=8$, independent of $r$.
Because a full rotation of the ellipse is obtained when $\Lambda_2(\tau)$
covers twice the unit circle, $\mathcal{A}_2=2 \pi$ and hence 
\[
   \geop = 16\pi \delta^2 + O(\delta^3),
\]
with contributions $32 \pi \delta^2$ and
$-16 \pi \delta^2$ from $\langle{\df \hat \phiform^{*}}\rangle$ and
$-\frac{1}{2}\langle{[\hat \phiform^{*},\hat \phiform^{*}]}\rangle$, 
respectively.
The discrepancy of $-2\pi$ when compared with \eqref{q:ellgeo2} results
from a  different definition of the angle $\theta$, which is measured
from an axis rotating with the ellipse in the calculation leading to
\eqref{q:ellgeo2} while it is measured from a fixed axis in \S6. 

%=============================================================================

% \bibliographystyle{elsart-num}
% \bibliography{all,dw}

\end{document}